# Implicit multistage two-derivative discontinuous Galerkin schemes for viscous conservation laws

Alexander Jaust, Jochen Schütz and David C. Seal

In this paper we apply implicit two-derivative multistage time integrators to conservation laws in one and two dimensions. The one dimensional solver discretizes space with the classical discontinuous Galerkin (DG) method, and the two dimensional solver uses a hybridized discontinuous Galerkin (HDG) spatial discretization for efficiency. We propose methods that permit us to construct implicit solvers using each of these spatial discretizations, wherein a chief difficulty is how to handle the higher derivatives in time. The end result is that the multiderivative time integrator allows us to obtain high-order accuracy in time while keeping the number of implicit stages at a minimum. We show numerical results validating and comparing methods.

## 1. Introduction

In this work, we focus on viscous conservation laws and present an implicit high-order time integration schemes for the discontinuous Galerkin (DG) method [11, 12, 13, 14, 15, 41]. One advantage of DG is that it is easy to increase spatial accuracy by locally increasing the polynomial degree, however this high-order accuracy is lost unless a high-order time integrator is applied to the result. For hyperbolic problems, explicit Runge-Kutta time integration schemes are often used because of their low computational costs, low dissipation, favorable stability regions and ease of implementation. However, the time step size $\Delta t$ is limited by Courant-Friedrichs-Lewy (CFL) stability constraints (which link $\Delta t$ to the mesh size and the speed of propagation of the physical system). In practical applications, the maximum allowable time step may be unacceptably small due to very small local mesh sizes. This becomes especially troublesome for problems with diffusion, where the maximum speed of propagation of information is infinite. These severe CFL restrictions can frequently be overcome by using implicit time integrators, which may permit for larger time steps at the expense of increased computational cost per time step. This is attributed to the fact that each time step requires one or more systems of (usually nonlinear) equations to be solved. Common examples of implicit time integrators include multistep methods such as the backward differentiation formulae (BDF) or multistage methods such as diagonally implicit Runge-Kutta (DIRK) methods [1].

One chief criticism of the discontinuous Galerkin method is the large number of degrees of freedom required to compute a solution. This is especially troublesome for the case of time dependent problems. As part of an effort to reduce the memory footprint for steady-state computations, the hybridized discontinuous Galerkin (HDG) method [10, 35, 42] was proposed as an alternative to more classical DG methods. More recently, this method has been applied to time-dependent problems [29, 34, 36, 35], which by construction, requires the use of implicit time integrators.

In this paper, we employ two-derivative methods as our time discretization. As the name suggests, these algorithms make use of more than one time derivative, and can be constructed to have a strong stability preserving (SSP) property [8, 37]. The principle can be best explained by starting with Taylor methods, which are a subclass of the so-called multistage multiderivative methods [23, 43]. If $w(x,t)$ is a function whose values are only known at time $t = t^n$, then approximations at some time $t > t^n$ can be obtained from a Taylor expansion in time:

$$w(x,t) = w(x,t^n) + (t-t^n)w_t(x,t^n) + \frac{(t-t^n)^2}{2}w_{tt}(x,t^n) + \cdots . \tag{1}$$

It is necessary to approximate additional time derivatives with this type of discretization. In the context of viscous conservation laws (that we define in Eqn. (2)), we replace temporal derivatives by spatial derivatives using the Cauchy-Kovalevskaya (CK) procedure.

In the context of numerical methods for partial differential equations, the use of Taylor series (in time) to discretize the PDE is often attributed to Lax and Wendroff [33]. There, the authors write down a second-order accurate Taylor series, and then appeal to the PDE to convert temporal derivatives to spatial derivatives. This approach is used for the so-called Taylor-Galerkin methods from the 1980's [17, 16], and the original ENO scheme of Harten et al. [25] uses the same procedure. In [18, 19], Dumbser and Munz construct discontinuous Galerkin schemes with arbitrary order of accuracy in space and time based on Arbitrary DERivative (ADER) schemes. They also present an approach to evaluate the Cauchy-Kovalevskaya procedure efficiently in [19] based on the work of Dyson [20] that relies on the application of the Leibniz rule. In addition, Qiu et al. [39] present an approach to couple Lax-Wendroff and discontinuous Galerkin methods (LWDG) based upon direct differentiation of the basis functions to define higher-order derivatives of the solution. They obtain a high-order, explicit one step method for hyperbolic problems that is up to third order accurate in time. They show that for their setting, the LWDG scheme is more efficient than a Runge-Kutta discontinuous Galerkin (RKDG). In [38] the approach is extended to 1D convection-diffusion equations based on the local discontinuous Galerkin (LDG) method. Furthermore, the behavior of the method coupled to different numerical fluxes is studied. Additionally, Taylor discretizations are investigated for finite difference weighted ENO methods in [7, 32, 40, 44].

Seal et. al. [43] are the first to extend the Lax-Wendroff type of approach to explicit multiderivative Runge-Kutta methods with DG and WENO spatial discretizations for hyperbolic conservation laws in a single dimension. They develop a framework for two-derivative Runge-Kutta methods that can be easily extended to incorporate additional stages or derivatives. In addition, Tsai et. al. [47] apply explicit and implicit two-derivative Runge-Kutta methods to PDEs with high-order finite-difference methods for spatial discretization.

Finally, we point out that the applications of the Cauchy-Kovalevskaya described thus far are used for temporal evolution of the solution, but the same procedure can also be used to define high-order boundary values by going the other way around. That is, in place of using the PDE to replace time derivatives with spatial derivatives, it is possible to define spatial derivatives (e.g., normal derivatives of the solution along the boundary of a domain) from time derivatives of the solution. This is done in the so-called inverse Lax-Wendroff (ILW) methods [27, 45, 46], as well as other related works [3, 26]. Our focus is not on this application, but on making use of the CK procedure to define the temporal evolution for the solution.

In this paper, we develop a strategy to apply implicit multistep two-derivative methods to convection-diffusion type equations in 1D using DG. That approach is then extended to first-order PDEs in 2D where we employ HDG for increased efficiency. The approach shows some similarities to the one in [47], but we use DG instead of finite differences for our spatial discretization. As the time derivatives that arise from two-derivative time integrators are replaced by spatial derivatives, an accurate way to represent them is needed. This could be done by differentiating the polynomial representation of the solution in the DG setting. This approach is justified in [43], which explains that the derivatives are multiplied by additional powers of $\Delta t$ that scale like $\Delta t = \mathcal{O}(\Delta x)$. In general, the time step restriction for implicit time integration is less severe. Therefore, we employ the LDG approach to accurately represent the additional derivatives, which has the additional benefit of potentially recovering superconvergence properties [21]. We refer the interested reader to [48], where the application of LDG to PDEs of higher order is discussed extensively. In this work, we also show how to extend our approach to efficiently solve for two dimensional convection problems with the hybridized DG method.

The remainder of this paper is structured as follows. In Section 2, we introduce a 1D nonlinear viscous conservation law that serves as a model equation. Then, in Section 3, we briefly describe the two-derivative multistage time integrators that are used in this work. This also explains the appearance of higher order spatial derivatives that are not directly present in the underlying PDE. Afterwards, we discretize the model equation in time and space using the LDG approach (c.f. Section 4) and verify the method using linear and nonlinear PDEs in Section 5. The two-derivative time integration is then extended to first order PDEs in two dimensions in Section 6. The resulting equation is discretized using an HDG method that significantly reduces the size of the globally coupled system. Finally, we verify the approach using the linear advection and nonlinear Euler equations in Section 7.

## 2. Underlying equation

In this work, we begin with the scalar nonlinear viscous conservation law

$$\begin{aligned} w_t + f(w)_x &= \varepsilon w_{xx} \quad \forall (x,t) \in \Omega \times \mathbb{R}^+ \\ w(x,0) &= w_0(x) \quad \forall x \in \Omega \end{aligned} \tag{2}$$

with $\varepsilon \geq 0$ given on a domain $\Omega \subset \mathbb{R}$ equipped with periodic boundary conditions. The method to be developed relies - similar to a Lax-Wendroff procedure [33] - on the use of the second temporal derivative $w_{tt}$, expressed in terms of spatial derivatives. For the underlying problem, we state the following lemma:

**Lemma 1.** *Let $w \in C^4(\Omega \times \mathbb{R}^+)$. Then, the second temporal derivative can be expressed as*

$$w_{tt} = \left(f'(w)f(w)_x - \varepsilon f'(w)w_{xx}\right)_x + \varepsilon\left(-f(w)_x + \varepsilon w_{xx}\right)_{xx} =: \mathcal{R}^2(w). \tag{3}$$

*Proof.* Obviously, there holds

$$w_t = -f(w)_x + \varepsilon w_{xx} =: \mathcal{R}(w)$$

and consequently,

$$\begin{aligned} w_{tt} &= (-f(w)_x + \varepsilon w_{xx})_t = (-f(w)_t)_x + \varepsilon (w_t)_{xx} = (-f'(w)w_t)_x + \varepsilon (w_t)_{xx} \\ &= \left(f'(w)f(w)_x - \varepsilon f'(w)w_{xx}\right)_x + \varepsilon\left(-f(w)_x + \varepsilon w_{xx}\right)_{xx} =: \mathcal{R}^2(w). \end{aligned}$$

□ □

**Remark 1** (Limiting cases). *The term for $w_{tt}$ simplifies significantly in some limiting cases:*

1. *If $f$ is linear, i.e., $f(w) = cw$, then*

$$w_{tt} = c^2 w_{xx} - 2\varepsilon c w_{xxx} + \varepsilon^2 w_{xxxx}.$$

2. *This also means that for $f \equiv 0$, $w_{tt} = \varepsilon^2 w_{xxxx}$.*
3. *If $\varepsilon \equiv 0$, there $w_{tt} = (f'(w)f(w)_x)_x$.*

*Note that the viscous and convective terms influence each other mutually, i.e., one obtains cross-terms that need to be dealt with.*

## 3. Time integration

In this section, we shortly review multiderivative time integrators as far as it is of importance for this work. Assume that the ordinary differential equation

$$y'(t) = g(y(t))$$

is given for a smooth function $g$. (As is customary, the prime $'$ denotes differentiation with respect to $t$.) Classical approaches (e.g., multistage Runge-Kutta, or linear multistep Adams methods) to the numerical approximation of these equations [22, 24] only use $g$ itself. A multiderivative method, on the other hand, takes knowledge of higher derivatives of $y$ into consideration. As an example, the second derivative $y''$ is given by

$$y''(t) = (g \circ y)'(t) := \partial_y g(y) \cdot y'(t) = \partial_y g(y) \cdot g(y(t)), \tag{4}$$

which can be readily computed for a system of ODEs using symbolic differentiation software. In this publication, we assume that $0 \le t \le T$, and that this temporal interval is uniformly subdivided into $0 = t^0 < t^1 < \ldots < t^N = T$ with spacing $\Delta t$. We note that none of the algorithms presented in this work depend on a uniform time step size; this choice is simply made for the ease of presentation. As is customary, $y^n$ denotes an approximation to $y$ at time $t = t^n$, $0 \le n \le N$.

The methods considered in this work are *implicit* two-point collocation methods that make use of multiple derivatives of the solution that take the form

$$\sum_{j=0}^{m} \Delta t^j (\partial_t^j y)(y^{n+1}) P^{(m-j)}(0) = \sum_{j=0}^{m} \Delta t^j (\partial_t^j y)(y^n) P^{(m-j)}(1), \tag{5}$$

where $P(t) = \frac{t^k (t-1)^\ell}{(k+\ell)!}$ and $\partial_t^j$ is the $j-$th temporal derivative of $y$. These methods can be found by fitting a Hermite-Birkhoff interpolation that matches a total of $k$ derivatives of the solution at time $t = t^n$, and $\ell$ derivatives of the solution at time $t = t^{n+1}$, and then integrating the result. In practice, at least for ODEs, an appropriate extension of (4) defines higher derivatives of the solution. Each of these methods are of order $m = \ell + k$ (cf. II.12 in [22]). Because of the growing complexity of higher order derivatives (see also Eqn. (3)), we rely on schemes involving only two derivatives of the unknown.

**Remark 2** (Employed methods)**.**

- *In general, the two point two-derivative method in this work can be written in the form*

$$\begin{aligned} y^{n+1} = y^n &+ \Delta t \left( \alpha_1 g(y^n) + \alpha_2 g(y^{n+1}) \right) \\ &+ \Delta t^2 \left( \beta_1 \dot{g}(y^n) + \beta_2 \dot{g}(y^{n+1}) \right), \end{aligned} \tag{6}$$

  *where the coefficients $\alpha_i, \beta_i$, for $i = 1, 2$ are chosen to increase the order of accuracy or modify the region of absolute stability of the method.*

- *In this work, we make use of the third-order method with $k = 1$, $\ell = 2$, given by*

$$y^{n+1} = y^n + \frac{\Delta t}{3} \left( g(y^n) + 2g(y^{n+1}) \right) - \frac{\Delta t^2}{6} \dot{g}(y^{n+1}), \tag{7}$$

  *and the fourth-order scheme with $k = \ell = 2$ given by*

$$y^{n+1} = y^n + \frac{\Delta t}{2} \left( g(y^n) + g(y^{n+1}) \right) + \frac{\Delta t^2}{12} \left( \dot{g}(y^n) - \dot{g}(y^{n+1}) \right). \tag{8}$$

  *These are the same methods that are used to discretize the non-linear terms in [5].*

**Lemma 2** (Stability)**.** *The integrators* (7) *and* (8) *are third- and fourth-order accurate, respectively, and A-stable. The third-order method* (7) *is L-stable.*

The order of accuracy is given in [22]. Here, we include a proof that these two methods are $A$, and $L$-stable, respectively. A plot of the stability region for the third-order method is presented in Fig. 1.

*Proof.* We apply each method to the test equation where $y' = \lambda y$, where $\lambda \in \mathbb{C}$ is a complex number. Method (7) results in

$$y^{n+1} = y^n + \frac{\Delta t}{3}\left(\lambda y^n + 2\lambda g(y^{n+1})\right) - \frac{\Delta t^2}{6}\lambda^2 y^{n+1}, \tag{9}$$

and method (8) becomes

$$y^{n+1} = y^n + \frac{\Delta t}{2}\left(\lambda y^n + \lambda y^{n+1}\right) + \frac{\Delta t^2}{12}\left(\lambda^2 y^n - \lambda^2 y^{n+1}\right). \tag{10}$$

If we define $\mu := \lambda \Delta t$, then each method can be written as

$$y^{n+1} = h(\mu) y^n \tag{11}$$

where

$$h(\mu) = \frac{1 + \frac{\mu}{3}}{1 - \frac{2}{3}\mu + \frac{\mu^2}{6}} \tag{12}$$

for method (7), and

$$h(\mu) = \frac{1 + \frac{\mu}{2} + \frac{\mu^2}{12}}{1 - \frac{\mu}{2} + \frac{\mu^2}{12}} \tag{13}$$

for method (8).

For method (8), we observe that $|h(iy)| = 1$ for any $y \in \mathbb{R}$, and that $\lim_{\mu\to\infty} |h(\mu)| = 1$. Because this function has no poles in the left half plane $\mathbb{C}^-$, the maximum modulus theorem indicates that $|h(\mu)| < 1$ for all $\mu \in \mathbb{C}^-$.

For method (7), we have $\lim_{\mu\to\infty} h(\mu) = 0$, and therefore this method has stiff decay. To obtain $L$-stability, we likewise need only show that $|h(iy)| \leq 1$ for any $y \in \mathbb{R}$ because this function has no poles in the left half plane. Omitting details for brevity, it can be shown that

$$|h(iy)|^2 = \frac{4(y^2+9)}{y^4+4y^2+36} = \frac{4y^2+36}{y^4+4y^2+36} \leq \frac{4y^2+36}{4y^2+36} = 1, \tag{14}$$

in which case $|h(iy))|^2 \leq 1$ for any $y \in \mathbb{R}$.

□

□

**Remark 3.** *One final observation is that the stability polynomials for the third- and fourth-order methods are identical to the Padé approximants $R_{1,2}(\mu)$, $R_{2,2}(\mu)$, respectively for $e^\mu$.*

Applying (6) to equation (2) on a semi-discrete level yields the expression

$$w^{n+1} - \Delta t \alpha_2 w_t^{n+1} - \beta_2 \Delta t^2 w_{tt}^{n+1} = w^n + \Delta t \alpha_1 w_t^n + \beta_1 \Delta t^2 w_{tt}^n,$$

where $w_{tt}$ has to be replaced by the expression in (3). The term $w_{tt}$ contains spatial derivatives up to fourth order, so we have to discuss how to discretize them in a DG framework efficiently. In [21, 48], the authors show how to use the *local discontinuous Galerkin* (LDG) method to discretize the higher spatial derivatives in an *explicit* DG solver. Their work will be the basis for the algorithm to be presented in the sequel.

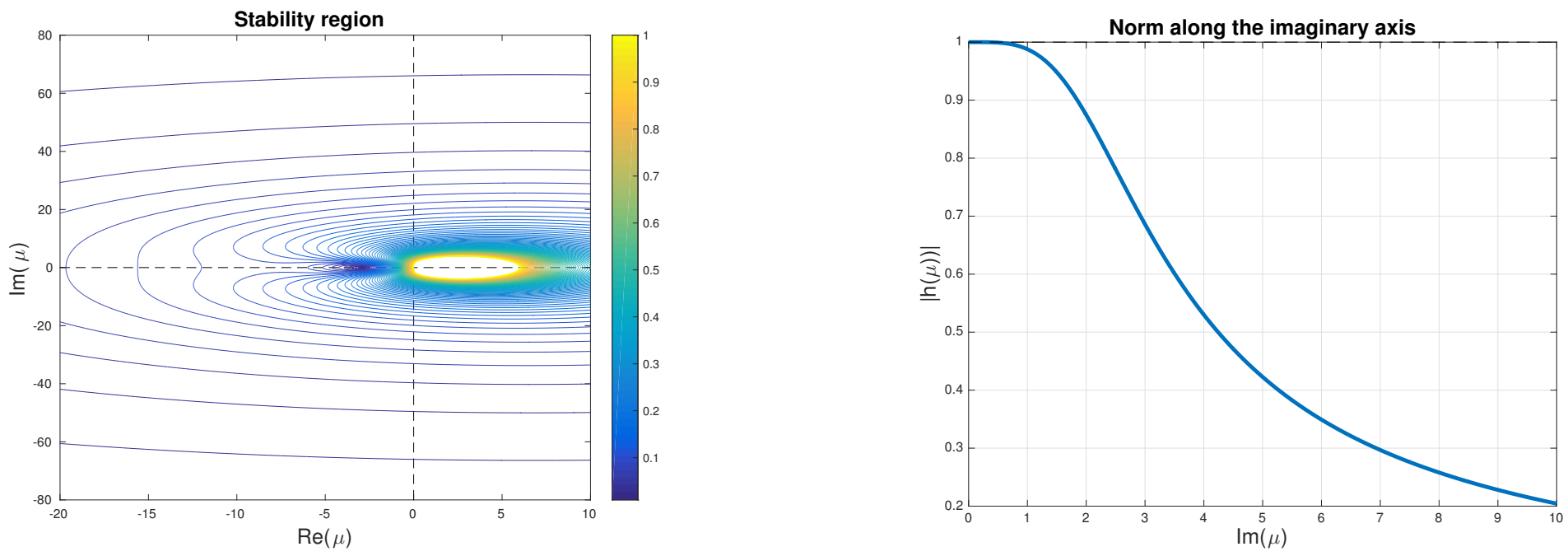

Figure 1: Plots of stability regions. Here, we plot the region of absolute stability for the $L$-stable third-order method defined in (7). In the left panel we plot $|h(\mu)|$ for various values of $\mu \in \mathbb{C}$, and in the right panel we plot $|h(iy)|$, where $h$ is defined in Eqn. (12), and $y$ is a real number. The stability region for the fourth order method is the left half plane $\mathbb{C}^-$, and therefore it is left out for brevity.

## 4. 1D: Spatial and temporal discretization

It is the aim of this publication to couple temporal discretization in (6) to the discontinuous Galerkin method. A semi-discretization of (2) is given by

$$\begin{aligned} w^{n+1} = w^n &+ \Delta t \left(\alpha_1 \mathcal{R}(w^n) + \alpha_2 \mathcal{R}(w^{n+1})\right) \\ &+ \Delta t^2 \left(\beta_1 \mathcal{R}^2(w^n) + \beta_2 \mathcal{R}^2(w^{n+1})\right), \end{aligned} \tag{15}$$

where $\mathcal{R}(w)$ and $\mathcal{R}^2(w)$ denote the expressions for $w_t$ and $w_{tt}$, respectively (c.f. Eqn. (3)), and $w^n$ denotes an approximation to $w$ at time $t^n$.

Before introducing the full spatial and temporal discretization, we start with some preliminaries. To introduce a finite element method, we begin by defining a triangulation of $\Omega$ into cells $\Omega_k$ such that they define a partition

$$\Omega = \bigcup_{k=1}^{N_e} \Omega_k$$

with a total of $N_e$ elements. For a given polynomial order $p$, we define the ansatz space $V_h$ to consist of cell-wise polynomials of order $p$ with no continuity restriction along the cell boundaries

$$V_h := \{q \in L^2(\Omega) \mid q \in \Pi^p(\Omega_k) \quad \forall k = 1, \ldots, N_e\}.$$

Again, it is possible to choose an adaptive $p$ that differs from cell to cell. We neglect this, for ease of exposition.

The method to be presented relies on the quantities

$$\sigma := w_x, \quad \tau := \sigma_x = w_{xx}, \quad \psi := \tau_x = w_{xxx}.$$

In the most straightforward way, these variables are discretized as

$$(\psi_h, \varphi_h)_{\Omega_k} + (\tau_h, (\varphi_h)_x)_{\Omega_k} - \langle \widehat{\tau}, \varphi_h n \rangle_{\partial\Omega_k} = 0 \quad \forall \varphi_h \in V_h, \tag{16}$$
$$(\tau_h, \varphi_h)_{\Omega_k} + (\sigma_h, (\varphi_h)_x)_{\Omega_k} - \langle \widehat{\sigma}, \varphi_h n \rangle_{\partial\Omega_k} = 0 \quad \forall \varphi_h \in V_h, \tag{17}$$
$$(\sigma_h, \varphi_h)_{\Omega_k} + (w_h, (\varphi_h)_x)_{\Omega_k} - \langle \widehat{w}, \varphi_h n \rangle_{\partial\Omega_k} = 0 \quad \forall \varphi_h \in V_h. \tag{18}$$

As is customary, we have defined the abbreviations

$$(f,g)_{\Omega_k} := \sum_{\Omega_k} \int_{\Omega_k} fg \mathrm{d}x, \quad \langle f,g\rangle_{\partial\Omega_k} := \sum_{\Omega_k} \int_{\partial\Omega_k} fg \mathrm{d}\sigma(x).$$

In one dimension, the last term can be simplified into function evaluations at two points. However, we prefer keeping the integral on the boundary to indicate the algorithm extends to multiple dimensions. The numerical fluxes $\widehat{\tau}$, $\widehat{\sigma}$ and $\widehat{w}$ have to be identified appropriately. One way to achieve a stable scheme is to choose 'upwinding' in an alternating fashion [48]. The corresponding fluxes read

$$\widehat{w} = w_h^+, \quad \widehat{\sigma} = \sigma_h^-, \quad \widehat{\tau} = \tau_h^+, \quad \widehat{\psi} = \psi_h^-, \tag{19}$$

where we stick to the convention that $w_h^-$ refers to the left value of $w_h$ at the interface, and $w_h^+$ refers to the right value.

We summarize these quantities in an auxiliary variable $\mathbf{x}_h \in V_h^4 =: \mathbf{X}_h$, given by

$$\mathbf{x}_h := (w_h, \sigma_h, \tau_h, \varphi_h). \tag{20}$$

This simplifies the defining equations (16)–(18) for $\sigma_h$, $\tau_h$ and $\varphi_h$ as

$$\mathcal{N}_{\text{aux}}(\mathbf{x}_h, \boldsymbol{\varphi}_h) = 0 \quad \forall \boldsymbol{\varphi}_h \in V_h^3. \tag{21}$$

For a convenient notation we use the following abbreviations

$$V_h^4 := V_h \times V_h \times V_h \times V_h, \; V_h^3 := V_h \times V_h \times V_h. \tag{22}$$

Note that this does not refer to the polynomial degree that is used.

**Remark 4** (Lifting operators)**.** *It is pointed out in [48] that it is possible to express the variable $\sigma$ in terms of $w$ via lifting operators, and subsequently $\tau$ and $\psi$ in terms of $w$ as well. This comes at the expense of computing lifting operators, i.e., to locally solve linear systems of equations in each cell (c.f. [2]).*

With these preliminaries, we now consider the semi-discretization (15) (see also (3)) once again. It is well-known how to spatially discretize $\mathcal{R}(w)$ using the DG method

$$(\mathcal{R}(w_h), \varphi_h)_{\Omega_k} \approx \mathcal{N}_{\mathcal{R}}(\mathbf{x}_h, \varphi_h) \tag{23}$$

$$:= (f(w_h) - \varepsilon\sigma_h, (\varphi_h)_x)_{\Omega_k} - \left\langle \widehat{f}(w_h^+, w_h^-) - \varepsilon\widehat{\sigma}, \varphi_h n \right\rangle_{\partial\Omega_k} \tag{24}$$

with the discretization of $\sigma$ and $\widehat{\sigma}$ as before. The numerical flux $\widehat{f}$ denotes a standard consistent and conservative Riemann solver. Details on the chosen flux are given in the numerical results section.

The discretization of $\mathcal{R}^2(w)$ as given in equation (3) is less straightforward. In particular, both the occurring higher derivatives and the nonlinearity of $f$ pose severe problems. Based on the definition of $\psi_h$, $\tau_h$ and $\sigma_h$ earlier, we propose the following discretization:

$$\left(\mathcal{R}^2(w_h), \varphi_h\right)_{\Omega_k} \approx \mathcal{N}_{\mathcal{R}^2}(\mathbf{x}_h, \varphi_h) := \tag{25}$$

$$- \left(f'(w_h)^2\sigma_h - \varepsilon f'(w_h)\tau_h, (\varphi_h)_x\right)_{\Omega_k} + \left\langle f'(\widehat{w}_h)^2\widehat{\sigma}_h - \varepsilon f'(\widehat{w}_h)\widehat{\tau}_h, \varphi_h n\right\rangle_{\partial\Omega_k} \tag{26}$$

$$+ \left(\varepsilon\mathcal{D}_h^2 f, (\varphi_h)_x\right)_{\Omega_k} - \left\langle \widehat{\varepsilon\mathcal{D}_h^2 f}, (\varphi_h)_x \right\rangle_{\partial\Omega_k} - \left(\varepsilon^2\psi_h, \varphi_h\right)_{\Omega_k} + \left\langle \varepsilon^2\widehat{\psi}_h, \varphi_h n \right\rangle_{\partial\Omega_k}. \tag{27}$$

Again, the fluxes $\widehat{w}$, $\widehat{\sigma}$, $\widehat{\tau}$ and $\widehat{\psi}$ are the LDG fluxes with alternating evaluation, see (19). $\mathcal{D}_h^2 f$ denotes an approximation to $f(w)_{xx}$, see Remark 5.

**Remark 5** (Discretization of $\mathcal{D}_h^2 f$)**.** *The suitable discretization of $\mathcal{D}_h^2 f \approx f(w)_{xx}$ depends on the choice of the convective flux $f$. We show two prototypical examples:*

1. *Linear equation, i.e., $f(w) = cw$. In this case, $f(w)_{xx} = cw_{xx}$, and a suitable choice is*

$$\mathcal{D}_h^2 f := c\tau, \qquad \widehat{\mathcal{D}_h^2 f} := c\widehat{\tau}.$$

2. *Burgers' equation, i.e., $f(w) = \frac{1}{2}w^2$. In this case, $f(w)_{xx} = w_x^2 + ww_{xx}$. As all occurring derivatives are known explicitly in the algorithm, a suitable approximation is*

$$\mathcal{D}_h^2 f := \sigma^2 + w\tau, \qquad \widehat{\mathcal{D}_h^2 f} := \widehat{\sigma}^2 + \widehat{w}\widehat{\tau}.$$

3. *A similar procedure as with Burgers' equation is possible with any flux function - also for Euler's equation. However, the result might become increasingly complex.*

Ultimately, this leads to the formulation of the full algorithm, summarized in the following definition:

**Definition 1** (Numerical method)**.** *Let $\boldsymbol{\varphi}_h = (\boldsymbol{\varphi}_h^{(1)}, \varphi_h^{(2)}) \in \mathbf{X}_h$ with $\boldsymbol{\varphi}_h^{(1)} \in V_h^3$ and $\varphi_h^{(2)} \in V_h$. Furthermore, let the semi-linear form $\mathcal{N}$ be given by*

$$\mathcal{N}(\mathbf{x}_h, \boldsymbol{\varphi}_h) := \begin{pmatrix} \mathcal{N}_{aux}(\mathbf{x}_h, \boldsymbol{\varphi}_h^{(1)}) \\ \mathcal{N}_{eq}(\mathbf{x}_h, \varphi_h^{(2)}) \end{pmatrix}$$

*where $\mathcal{N}_{eq}(\mathbf{x}_h, \varphi_h^{(2)})$ is given by*

$$\begin{aligned} \mathcal{N}_{eq}(\mathbf{x}_h, \varphi_h^{(2)}) := {} & \alpha_1 \mathcal{N}_{\mathcal{R}}(\mathbf{x}_h^n, \varphi_h^{(2)}) + \alpha_2 \mathcal{N}_{\mathcal{R}}(\mathbf{x}_h^{n+1}, \varphi_h^{(2)}) \\ & + \Delta t \left( \beta_1 \mathcal{N}_{\mathcal{R}^2}(\mathbf{x}_h^n, \varphi_h^{(2)}) + \beta_2 \mathcal{N}_{\mathcal{R}^2}(\mathbf{x}_h^{n+1}, \varphi_h^{(2)}) \right). \end{aligned}$$

*The coefficients $\alpha_i$ and $\beta_i$ are the same as in Remark 2 and are chosen to modify order of accuracy or stability of the time integrator. The approximate solution $\mathbf{x}_h^{n+1} = (w_h^{n+1}, \sigma_h^{n+1}, \tau_h^{n+1}, \psi_h^{n+1}) \in \mathbf{X}_h$ is given as the solution to the problem*

$$\begin{pmatrix} \mathbf{0} \\ \frac{1}{\Delta t} \left( w_h^{n+1} - w_h^n, \varphi_h^{(2)} \right)_{\Omega_k} \end{pmatrix} = \mathcal{N}(\mathbf{x}_h, \boldsymbol{\varphi}_h) \quad \forall \boldsymbol{\varphi}_h \in \mathbf{X}_h.$$

*Note that the first component (which is indeed vector-valued in $\mathbb{R}^3$, i.e., $\mathbf{0} \in \mathbb{R}^3$) stems from* (16)–(18)*, while the second component is the discretized version of equation* (15)*.*

The following lemma is a straightforward consequence of both the order of accuracy of the ODE integrator and the consistency of the underlying DG schemes:

**Lemma 3** (Consistency in time)**.** *The algorithm is consistent with the order of the temporal integration scheme chosen in* (15)*, i.e., there holds:*

$$\begin{pmatrix} \mathbf{0} \\ \frac{1}{\Delta t} \left( w(\cdot, t^{n+1}) - w(\cdot, t^n), \varphi_h \right)_{\Omega_k} \end{pmatrix} - \mathcal{N}(\mathbf{x}, \boldsymbol{\varphi}_h) = \mathcal{O}(\Delta t^q)$$

*where $q = 3$ for integrator* (7) *and $q = 4$ for integrator* (8)*.*

**Lemma 4** (Conservation)**.** *The algorithm is both locally and globally conservative if $\widehat{\mathcal{D}_h^2 f}$ is conservative.*

*Proof.* Testing with a piecewise constant test function yields that the integral of $w_h^{n+1}$ only depends on the fluxes over the boundaries. This yields local conservation. Noting that the fluxes are conservative and testing against a constant function yields that the algorithm is globally conservative. □ □

## 5. Numerical results: 1D examples

In this section, we present numerical results for the newly developed scheme. In each case, we demonstrate that the optimal order of convergence is met.

In all our computations, we use periodic boundary conditions on the unit interval $\Omega := [0, 1]$, and compute until a final time of $T = 0.5$. For the cases involving linear convection, we choose the upwind numerical flux

$$\widehat{f}(w^+, w^-) := cw^-, \quad c > 0,$$

whereas for Burgers equation, we use a local Lax-Friedrichs flux. The domain $\Omega$ is subdivided into equally spaced intervals with spacing $h$. As an error measure, we compute the $L^2$-error at time $T$, that is, we define the error as

$$e_h := \|w(\cdot, T) - w_h(\cdot, T)\|_{L^2(\Omega)},$$

where $w$ is the exact and $w_h$ the approximate solution to the underlying problem.

### 5.1. Heat equation

The first equation to be considered is the pure heat equation

$$w_t = \varepsilon w_{xx} \quad \forall (x, t) \in \Omega \times (0, T)$$

with initial conditions $w_0(x) = \sin(2\pi x)$, and $\varepsilon = 0.1$.

Numerical results for different values of the polynomial order $p$ of the ansatz space are shown in Fig. 2 for the third-order integrator (7) (left) and the fourth-order integrator (8) (right). The expected order of accuracy of $\max\{p + 1, 3\}$ and $\max\{p + 1, 4\}$, respectively, is achieved. The time step is set to $\Delta t = \Delta x$. Experiments with other ratios of $\frac{\Delta t}{\Delta x}$ introduce no stability problems, which is independent on the choice of $\varepsilon$. Thus, we conjecture that the algorithm is uniformly stable for this simple 1D test case without transport.

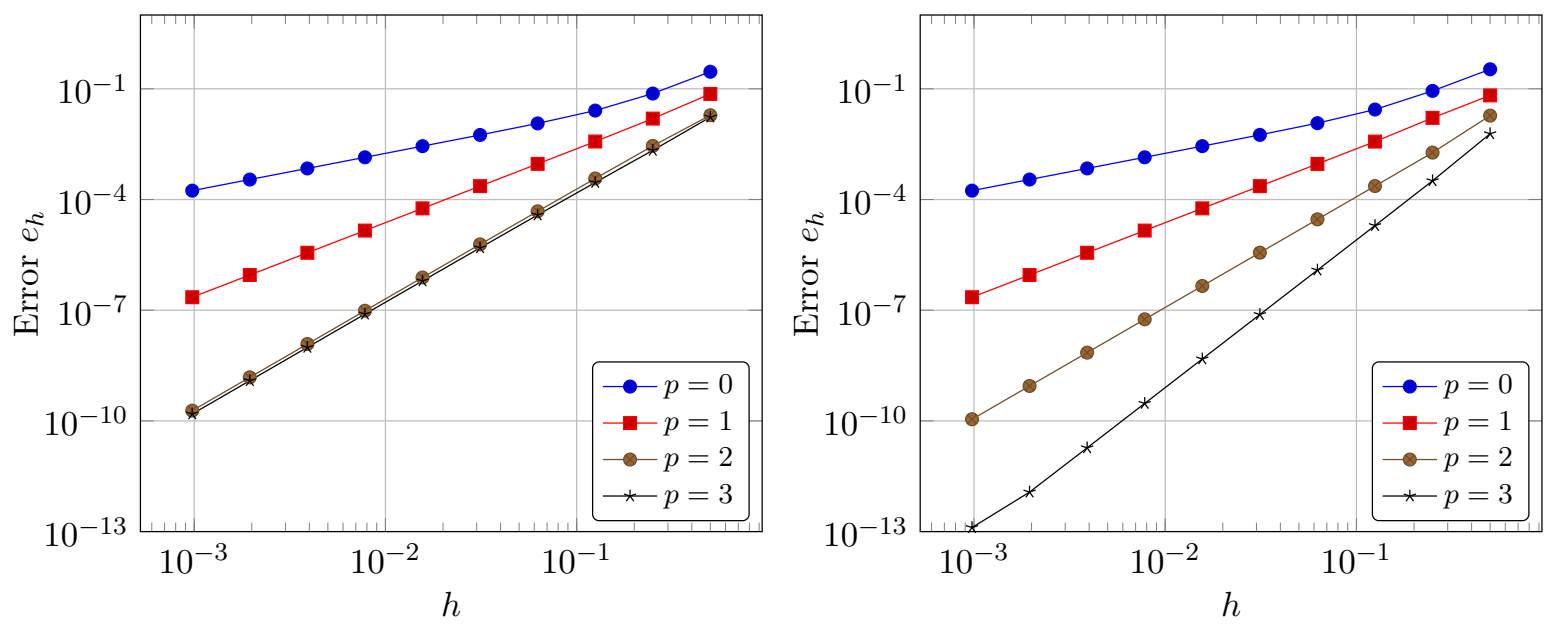

Figure 2: Numerical results for the heat equation. In both computations, we choose the ratio $\frac{\Delta t}{\Delta x}$ to be one. Temporal integration is performed via the third-order accurate integrator (7) (left) and the fourth-order integrator (8) (right). Tabulated results are given in Tab. 1.

### 5.2. Convection equation

Next, we test the algorithm on the pure convection equation

$$w_t + cw_x = 0 \quad \forall (x, t) \in \Omega \times (0, T),$$

again with initial conditions $w_0(x) = \sin(2\pi x)$, and constant $c = 1$. Numerical results are displayed in Fig. 3, again for the third-order (left) and the fourth-order (right) temporal integrator. The CFL number for this example is one for the third-order integrator, and only 0.1 for the fourth-order integrator. The reason for this choice is that we find stability constraints with the fourth-order integrator. Our experience with other time integrators has lead to this in the past, and we suspect that it is most likely due to the loss of $L-$stability in the fourth-order solver. Numerical experiments indicate that the third-order integrator is uniformly stable. We note that obviously, such a severe CFL restriction is not a desired feature of an implicit scheme, and a detailed investigation into how to fix the fourth-order scheme is the subject of future work.

For this example, we also ran experiments for a longer time $T$. For brevity, we do not report the results here, but simply state that the methods behave as expected: optimal orders of accuracy are met.

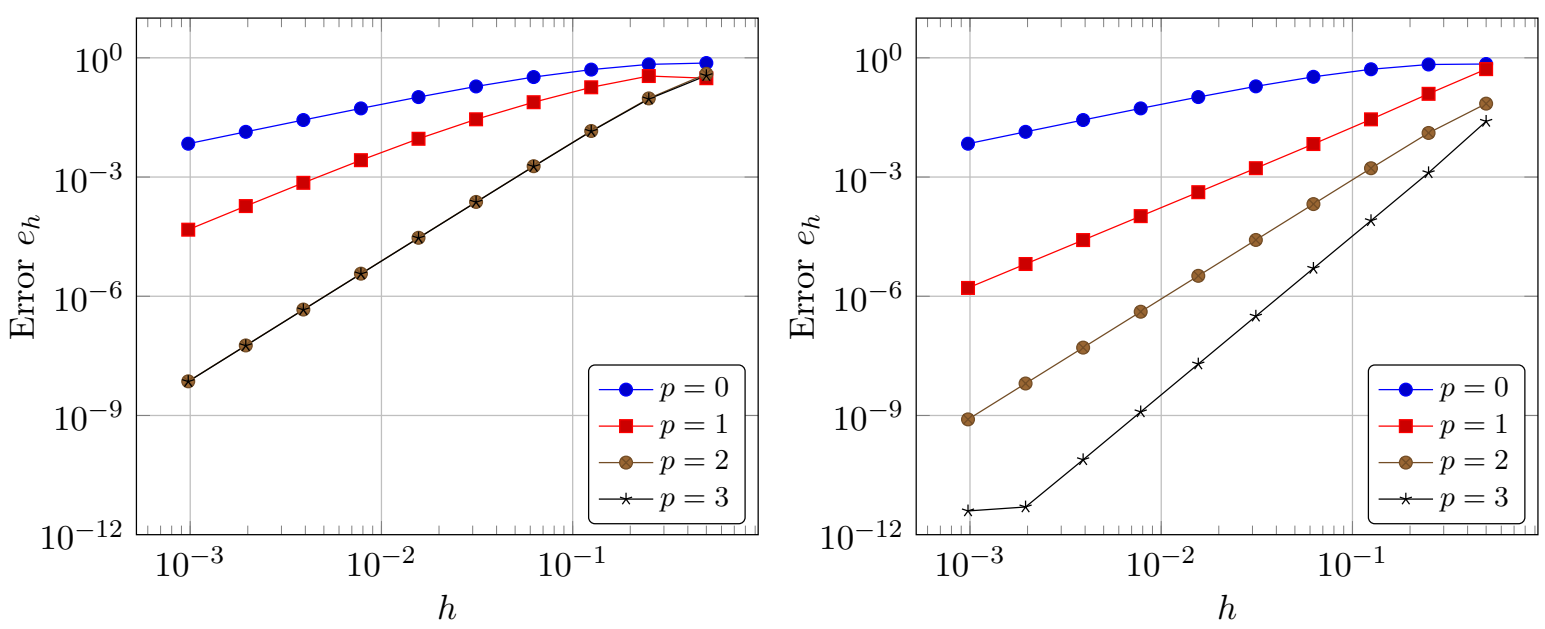

Figure 3: Numerical results for the convection equation. Temporal integration is done via the third-order integrator (7) (left) and the fourth-order integrator (8) (right). The ratio $\frac{\Delta t}{\Delta x}$ is set to be 1.0 (left) and 0.1 (right). The third-order integrator seems to be uniformly stable, yet the fourth-order integrator is not, which is why we find it necessary to reduce the CFL number. Tabulated results are given in Tab. 2.

### 5.3. Convection-diffusion equation

The final linear single-dimensional test case is the convection-diffusion equation

$$w_t + cw_x = \varepsilon w_{xx} \quad \forall (x,t) \in \Omega \times (0,T),$$

with values $c = 1$ and $\varepsilon = 0.1$. This exercises the ability of the scheme to correctly account for the additional coupling terms that arise in the discretization of $w_{tt}$. We present two examples: a) an example with a smooth initial profile, and b) a problem with discontinuous initial conditions.

#### 5.3.1. Convection-diffusion: Smooth initial conditions.

The initial conditions for this test problem are given by $w_0(x) = \sin(2\pi x)$, and the ratio $\frac{\Delta t}{\Delta x}$ is chosen to be one for all cases. Numerical results are presented in Figs. 4 and 5. In Fig. 4, we demonstrate a convergence study for various values of the polynomial degree $p$. In Fig. 5, we compare errors versus time to solution against other competing time stepping methods. We choose two implicit DIRK methods to compare against because we have used these methods in the past and found that they are efficient [29]. Moreover, because they are implicit methods written in the same framework, we can make use of identical data structures in order to provide fair timing results. The fourth-order method we use can be found in Hairer and Wanner's text [24], and the third-order method is attributed to Cash [4]. To keep the plot from becoming too complicated, we only show results for $p = 3$, however, all other results look similarly. One can observe that both the error levels and also the timing results are more or less the same. In particular, for large mesh sizes, it seems that the multiderivative method has a small advantage over DIRK methods.

For the third-order integrator, we do not observe any stability issues, while obviously, from our experiences with the convection equation, the fourth-order scheme tends to be only conditionally stable.

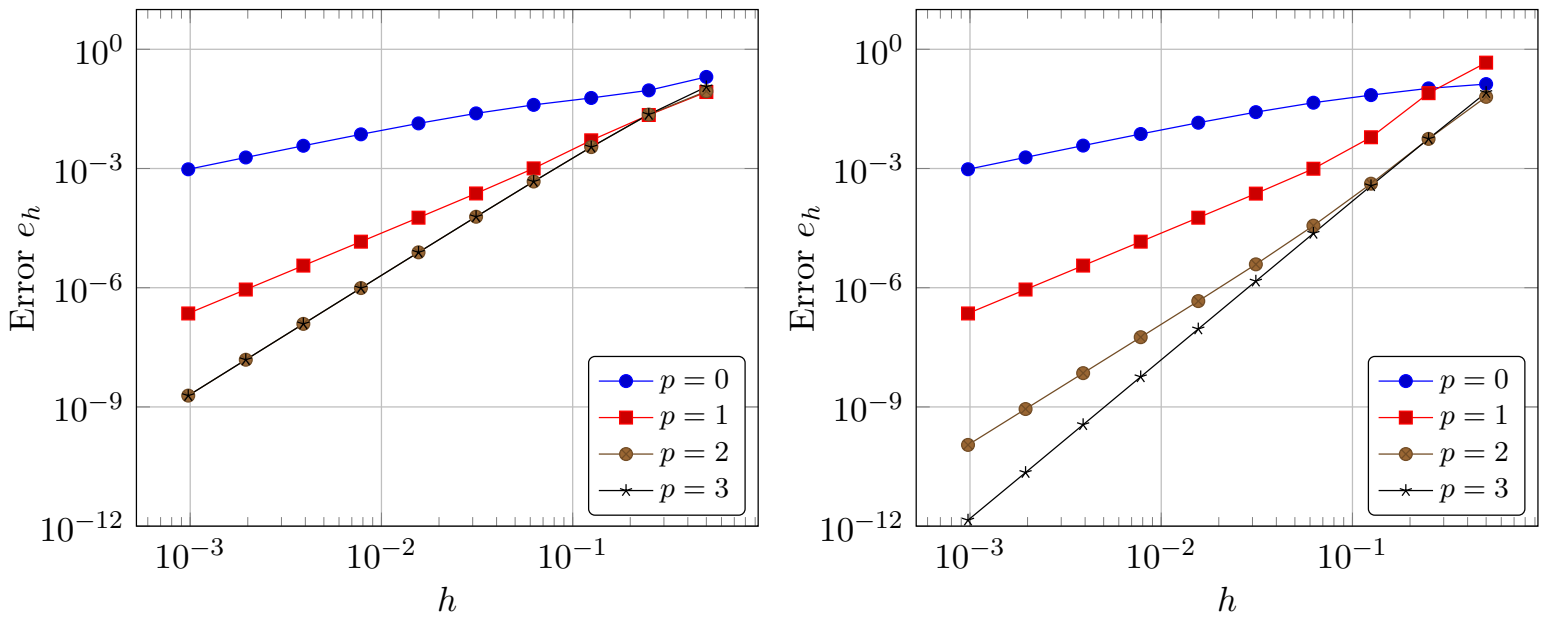

Figure 4: Numerical results for the convection-diffusion equation. In both computations, we choose the ratio $\frac{\Delta t}{\Delta x}$ to be one. Temporal integration is done via the third-order integrator (7) (left) and the fourth-order integrator (8) (right). Tabulated results are given in Tab. 3.

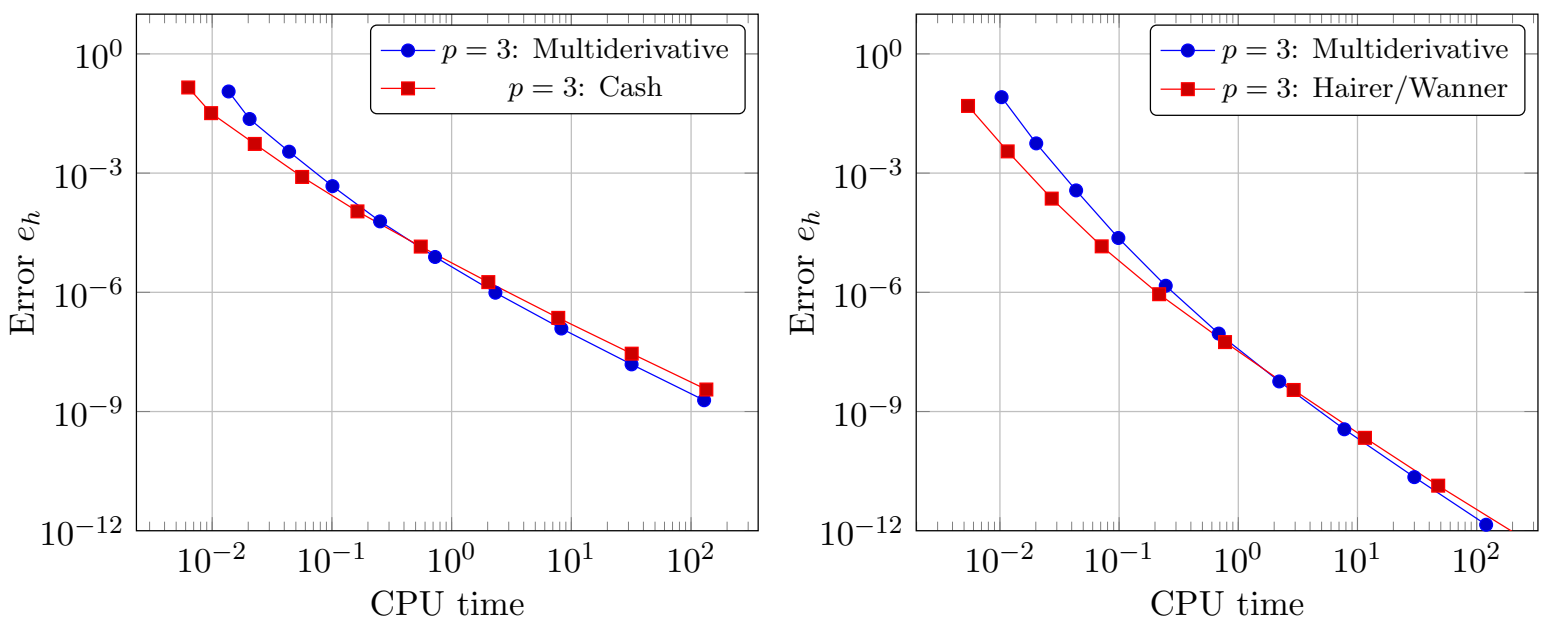

Figure 5: Timing results for the convection-diffusion equation. The settings are the same as in Fig. 4, however this time, we plot error versus computational time. Third order plots (left) use Cash's DIRK method, fourth-order plots use Hairer/Wanner's DIRK method. The plots indicate that both error levels and timings are very comparable. Tabulated results are given in Tab. 4.

### 5.3.2. Convection-diffusion: Discontinuous initial conditions.

In this section, we consider discontinuous initial data

$$w_0(x) = \mathcal{H}(\sin(2\pi(x-0.3)))e^{\sin(2\pi x)},$$

where $\mathcal{H}$ denotes the Heaviside function.

Because the initial conditions are not smooth, there are at least two ways of defining the initial conditions for the auxiliary variables $\mathbf{x}_h^0$ in Eqn. (20) that require the spatial derivatives of the initial conditions. One way is to fill it with the given initial conditions, i.e., set $\sigma_h^0 = \Pi_{V_h} w_0' \ldots$ where $\Pi_{V_h}$ denotes the $L^2-$projection onto $V_h$. We use this choice for the computations in the previous sections, because the initial datum is smooth. However, this problem has non-smooth initial conditions, and therefore $w_0'$ does not exist. (At least not in the classical sense, as it is a Dirac measure.) For this case, one alternative is to compute $\sigma_h^0$, $\tau_h^0$ and $\psi_h^0$ as a solution to $\mathcal{N}_{\text{aux}}(\mathbf{x}_h^0, \boldsymbol{\varphi}_h) = 0$ (see also (16)-(18) for the defining equations) for all $\boldsymbol{\varphi}_h \in V_h^3$ for a given $w_h = \Pi_{V_h} w_0$. This is the ansatz we pursue in this section.

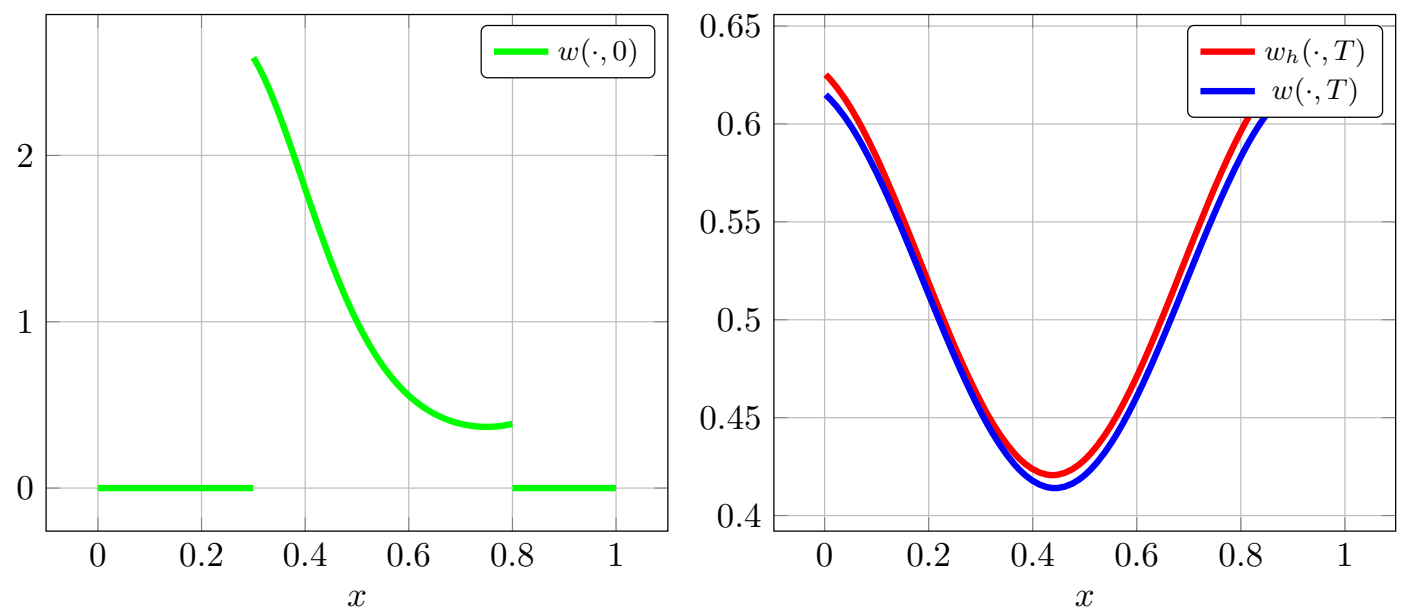

Figure 6: Approximate versus exact solution at time $T = 0.5$ (right) to the convection-diffusion equation with discontinuous initial datum (left). Results are computed with $p = 2$, $h = \frac{1}{16}$ and the third-order integrator (7).

In Fig. 6, we show an approximate solution at time $T = 0.5$, which uses a spatial width of $h = \frac{1}{16}$, quadratic (third-order) polynomials, CFL number of $\frac{\Delta t}{\Delta x} = 0.5$ and the third-order ODE integrator. We observe a strong agreement between the exact and approximate solution.

### 5.4. Viscous Burgers equation

Our final single dimensional numerical result is the nonlinear Burgers equation

$$w_t + f(w)_x = \varepsilon w_{xx}$$

with $f(w) = 0.5w^2$ and $\varepsilon = 0.1$. Equipped with initial conditions $w_0(x) = \sin(2\pi x)$, this test case has a smooth solution $w$ for all times $T$. As before, we show convergence results in Fig. 7. The exact solution is computed using the Cole-Hopf transformation [28]. No stability issues are observed in the computations, and the plots show perfect order of convergence. The results are similar to those of the convection-diffusion equation, which is mainly because diffusion is dominant in this test case. For implicit methods, this is probably the most relevant case, as for purely hyperbolic problems, explicit methods are often times the preferred method of choice given the finite speed of propagation of information. Furthermore, we do not implement any limiters but a final algorithm should have a suitable way of stabilizing discontinuities in the case of a convection dominated problem. This is one topic of future research, where one option is to introduce artificial viscosity into the time stepping. This is appealing because we already have access to these higher derivatives.

## 6. 2D: Extensions to multiple dimensions

In this section, we describe the extension of implicit two-derivative multistage methods to hyperbolic systems of first-order PDEs

$$\begin{aligned} \vec{w}_t + \nabla \cdot \vec{f}(\vec{w}) &= 0 \qquad \forall (\vec{x}, t) \in \Omega \times \mathbb{R}^+ \\ \vec{w}(\vec{x}, 0) &= \vec{w}_0(\vec{x}) \quad \forall \vec{x} \in \Omega \end{aligned} \tag{28}$$

on a domain $\Omega \subset \mathbb{R}^2$ with appropriate boundary conditions. In general, the unknown $\vec{w}$ is a function of space and time $\vec{w} := \vec{w}(\vec{x}, t)$, but we may drop $\vec{x}$ and $t$ for a more compact notation. We allow $\vec{f}(\vec{w})$ to be a general, possibly nonlinear, flux. Note, that the flux $\vec{f}(\vec{w})$ (as well as $\vec{w}$) is a vector in $\mathbb{R}^m$ for a system with a total of $m$ equations. In order to limit the complexity we stay with first-order PDEs for the 2D case. Due to the additional spatial dimension, the total number of derivatives also doubles which would make the assembly of the matrices required for second or higher order PDEs tedious, especially in the case of a system of equations. Nevertheless, the two-derivative methods are still applicable to higher order PDEs in 2D using the techniques from Section 4.

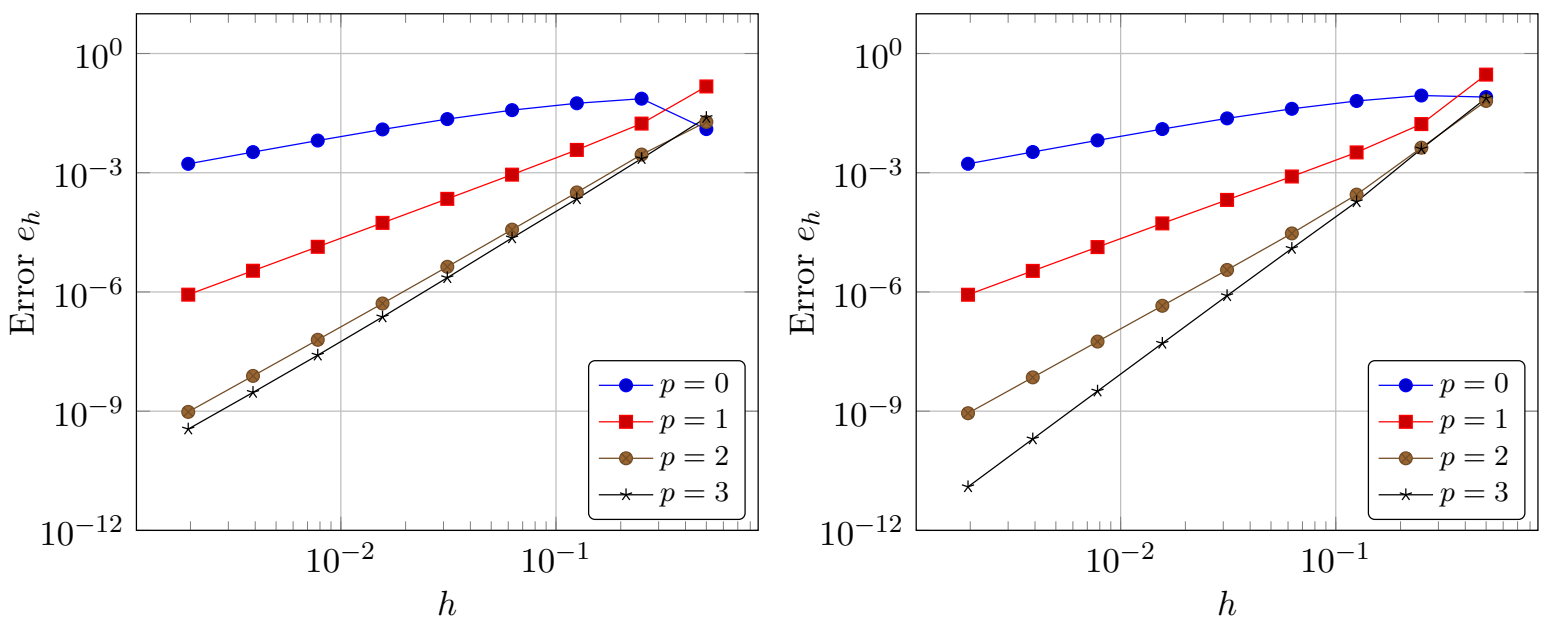

Figure 7: Numerical results for Burgers equation. In both computations, we choose the ratio $\frac{\Delta t}{\Delta x}$ to be one. Temporal integration is done via the third-order integrator (7) (left) and the fourth-order integrator (8) (right). Tabulated results are given in Tab. 5.

Here, we apply the same third and fourth-order two derivative methods defined in equations (7)-(8), but we discretize the resulting system using the *hybridized discontinuous Galerkin* [10, 35, 42] method. For Poisson's equation, this discretization is equivalent to (a variant of) the LDG method [10]. A preliminary investigation that couples two-derivative Runge-Kutta methods with HDG for a linear advection equation can be found in in [30]. On the one side, the HDG method really depends on efficient implicit time integrators due to the stiffness of the system of equations. Therefore, it is especially important to find new implicit time integrators that may be better than the currently employed BDF and DIRK methods. On the other side, the HDG method usually leads to a much smaller system of globally coupled equations than the LDG approach. This is especially beneficial in the current case where additional unknowns are introduced by the spatial derivatives from the two-derivative time discretization.

In this section, we only consider convection equations, treating diffusive parts is left for future work, as the treatment of higher-order derivatives in the HDG method is by far not standard. (See [6, 9] for the extension of HDG to higher-order derivatives.)

The representations of $\mathcal{R}(\vec{w})$ and $\mathcal{R}^2(\vec{w})$ only differ slightly from the one dimensional DG case. Here, we find that the first and second derivatives are given by

$$\begin{aligned} \vec{w}_t &= -\nabla \cdot \vec{f}(\vec{w}) =: \mathcal{R}(\vec{w}) \\ \vec{w}_{tt} &= \nabla \cdot (\vec{f}'(\vec{w}) \nabla \cdot \vec{f}(\vec{w})) =: \mathcal{R}^2(\vec{w}) \end{aligned} \tag{29}$$

which follows directly from Eqn. (3) by setting $\varepsilon = 0$.

We follow the discretization procedure of Section 4. The semi-discrete system again reads

$$\vec{w}^{n+1} = \vec{w}^n + \Delta t \left(\alpha_1 \mathcal{R}(\vec{w}^n) + \alpha_2 \mathcal{R}(\vec{w}^{n+1})\right) + \Delta t^2 \left(\beta_1 \mathcal{R}^2(\vec{w}^n) + \beta_2 \mathcal{R}^2(\vec{w}^{n+1})\right)$$

with $\mathcal{R}(\vec{w})$ and $\mathcal{R}^2(\vec{w})$ defined in Eqn. (29). For the spatial discretization we triangulate the domain such that

$$\Omega = \bigcup_{k=1}^{N_e} \Omega_k.$$

The hybridized DG method requires a description of the edges. We refer to edges of two intersecting elements and elements intersecting the domain boundary $\partial\Omega$ with $e_k$. The set of all edges is $\Gamma$ and its number of elements is $N_f := |\Gamma|$. This is needed to introduce a new hybrid unknown $\vec{\lambda} = \vec{w}_{|\Gamma}$ that represents the solution evaluated on the trace of each element. This allows us to reduce the size of the globally coupled system by using static condensation [10]. For the approximation of $\vec{\lambda}_h \approx \vec{\lambda}$ we need to introduce the ansatz space $M_h$ that consists of edge-wise polynomials of degree $p$ defined by

$$M_h := \{q \in L^2(\Gamma) \mid q_{|e_k} \in \Pi^p(e_k) \forall k = 1, \ldots, N_f,\ e_k \in \Gamma\}^m.$$

For the approximation of second order spatial derivatives introduced by a two-derivative time discretization, we again define an auxiliary variable through

$$\vec{\sigma} := \nabla \vec{w}.$$

The ansatz spaces for $\vec{\sigma}_h$ and $\vec{w}_h$ are the common spaces

$$\begin{aligned} H_h &:= \{q \in L^2(\Omega) \mid q_{|\Omega_k} \in \Pi^p(\Omega_k)\ \forall k = 1, \ldots, N_e\}^{2m}, \\ V_h &:= \{q \in L^2(\Omega) \mid q_{|\Omega_k} \in \Pi^p(\Omega_k)\ \forall k = 1, \ldots, N_e\}^{m}, \end{aligned}$$

that contain all polynomials of degree at most $p$. In order to condense notation, we define the vector of unknowns as

$$\mathbf{x}_h := (\vec{\sigma}_h, \vec{w}_h, \vec{\lambda}_h)$$

that stems from the ansatz space $\mathbf{X}_h := H_h \times V_h \times M_h$, and corresponding test functions $\boldsymbol{\varphi}_h = (\vec{\varphi}^{(1)}, \vec{\varphi}^{(2)}, \vec{\varphi}^{(3)}) \in \mathbf{X}_h$. Then, $\vec{\sigma}_h$ is approximated through

$$\left(\vec{\sigma}_h, \vec{\varphi}_h^{(1)}\right)_{\Omega_k} + \left(\vec{w}_h, \nabla \cdot (\vec{\varphi}_h^{(1)})\right)_{\Omega_k} - \left\langle \vec{\lambda}_h, \vec{\varphi}_h^{(1)} \cdot \vec{n} \right\rangle_{\partial\Omega_k} = 0 \quad \forall \vec{\varphi}_h^{(1)} \in H_h,$$

that is very similar to the approximation given previously, but we use the hybrid variable $\vec{\lambda}_h$ as the numerical flux $\widehat{\vec{w}} := \vec{\lambda}_h$. The equation is abbreviated by

$$\mathcal{N}_{\text{aux}}(\mathbf{x}_h, \vec{\varphi}_h^{(1)}) = 0 \quad \forall \vec{\varphi}_h^{(1)} \in H_h. \tag{30}$$

Finally, the discretization of $\mathcal{R}$ and $\mathcal{R}^2$ for this first-order PDE is given by

$$\left(\mathcal{R}(\vec{w}_h), \vec{\varphi}_h^{(2)}\right)_{\Omega_k} \approx \mathcal{N}_{\mathcal{R}}(\mathbf{x}_h, \vec{\varphi}_h^{(2)}) := \left(\vec{f}(\vec{w}_h), \nabla \vec{\varphi}_h^{(2)}\right)_{\Omega_k} - \left\langle \widehat{\vec{f}}, \vec{\varphi}_h^{(2),-} \vec{n} \right\rangle_{\partial\Omega_k}$$

and

$$\begin{aligned} \left(\mathcal{R}^2(\vec{w}_h), \vec{\varphi}_h^{(2)}\right)_{\Omega_k} &\approx \mathcal{N}_{\mathcal{R}^2}(\mathbf{x}_h, \vec{\varphi}_h^{(2)}) \\ &:= -\left(\mathcal{D}_h \vec{f}(\vec{w}_h, \vec{\sigma}_h), \nabla \vec{\varphi}_h^{(2)}\right)_{\Omega_k} + \left\langle \widehat{\mathcal{D}_h \vec{f}}, \vec{\varphi}_h^{(2),-} \vec{n} \right\rangle_{\partial\Omega_k}, \end{aligned}$$

where

$$\vec{f}'(\vec{w}) \nabla \cdot \vec{f}(\vec{w}) = \vec{f}'(\vec{w}) \vec{f}_i'(\vec{w}) \partial_{x_i} \vec{w} \approx \vec{f}'(\vec{w}_h) \vec{f}_i'(\vec{w}_h) \vec{\sigma}_{h,i} =: \mathcal{D}_h \vec{f}(\vec{w}_h, \vec{\sigma}_h).$$

Thus, $\vec{\sigma}_h$ is involved when the flux is evaluated. On each element interface, we insert numerical fluxes

$$\begin{aligned} \widehat{\vec{f}} &= \vec{f}(\vec{\lambda}_h) + \eta(\vec{w}_h^- - \vec{\lambda}_h)\vec{n} \\ \widehat{\mathcal{D}_h \vec{f}} &= \widehat{\mathcal{D}_h \vec{f}}(\vec{\lambda}_h, \vec{\sigma}_h^-) - \theta(\vec{w}_h^- - \vec{\lambda}_h)\vec{n}, \end{aligned}$$

that are modified Lax-Friedrichs (Rusanov) fluxes with $\eta$ and $\theta$ being stabilization parameters. Whenever $\theta$ is multiplied with a negative coefficient from the time integrator, we invert the sign. Note, that at this point the equations are only coupled through the hybrid variable $\vec{\lambda}_h$. An additional equation arises from the additional unknown $\vec{\lambda}_h$ through

$$\left\langle -\widehat{\vec{f}} + \widehat{\mathcal{D}_h \vec{f}}, \vec{\varphi}_h^{(3)} \right\rangle_{\partial\Omega_k} = 0 \quad \forall \vec{\varphi}_h^{(3)} \in M_h.$$

With these preliminaries in place, we are now ready to define the full hybridized DG method.

**Definition 2** (HDG method). *Let $\vec{\varphi}_h = \mathbf{X}_h$. Furthermore, let the semi-linear form $\mathcal{N}$ be given by*

$$\mathcal{N}(\mathbf{x}_h, \vec{\varphi}_h) := \begin{pmatrix} \mathcal{N}_{aux}(\mathbf{x}_h, \vec{\varphi}_h^{(1)}) \\ \mathcal{N}_{eq}(\mathbf{x}_h, \vec{\varphi}_h^{(2)}) \\ \mathcal{N}_{hyb}(\mathbf{x}_h, \vec{\varphi}_h^{(3)}) \end{pmatrix},$$

*with $\mathcal{N}_{aux}(\mathbf{x}_h, \vec{\varphi}_h^{(1)})$ as defined in equation (30), where $\mathcal{N}_{eq}(\mathbf{x}_h, \vec{\varphi}_h^{(2)})$ is given by*

$$\begin{aligned} \mathcal{N}_{eq}(\mathbf{x}_h, \vec{\varphi}_h^{(2)}) :=& \alpha_1 \mathcal{N}_{\mathcal{R}}(\mathbf{x}_h^n, \vec{\varphi}_h^{(2)}) + \alpha_2 \mathcal{N}_{\mathcal{R}}(\mathbf{x}_h^{n+1}, \vec{\varphi}_h^{(2)}) \\ &+ \Delta t \left( \beta_1 \mathcal{N}_{\mathcal{R}^2}(\mathbf{x}_h^n, \vec{\varphi}_h^{(2)}) + \beta_2 \mathcal{N}_{\mathcal{R}^2}(\mathbf{x}_h^{n+1}, \vec{\varphi}_h^{(2)}) \right), \end{aligned}$$

*and the hybrid term is given by*

$$\mathcal{N}_{hyb}(\mathbf{x}_h, \vec{\varphi}_h^{(3)}) := \left\langle [\![ -\alpha_2 \widehat{\vec{f}}^{\,n+1} + \beta_2 \widehat{\mathcal{D}_h \vec{f}}^{\,n+1} ]\!], \vec{\varphi}_h^{(3)} \right\rangle_{\partial\Gamma}.$$

*The brackets denote the jump operator*

$$[\![\vec{v}]\!] = \vec{v}^- \vec{n} - \vec{v}^+ \vec{n}$$

*with $v(\vec{x})^\pm$ being*

$$\vec{v}(\vec{x})^\pm = \lim_{\epsilon \to 0} \vec{v}(\vec{x} \pm \epsilon \vec{n}), \; \vec{x} \in \partial\Omega_k \tag{31}$$

*where $\vec{n}$ is the outward pointing normal. The approximate solution $\mathbf{x}_h^{n+1} = (\vec{\sigma}_h^{n+1}, \vec{w}_h^{n+1}, \vec{\lambda}_h^{n+1}) \in \mathbf{X}_h$ is given as the solution to the problem*

$$\begin{pmatrix} \mathbf{0} \\ \frac{1}{\Delta t} \left( \vec{w}_h^{n+1} - \vec{w}_h^n, \vec{\varphi}_h^{(2)} \right)_{\Omega_k} \\ \mathbf{0} \end{pmatrix} = \mathcal{N}(\mathbf{x}_h, \vec{\varphi}_h) \quad \forall \vec{\varphi}_h \in \mathbf{X}_h.$$

**Remark 6** (Number of unknowns). *All evaluations on elements only depend on local values of $\vec{w}_h$ and $\vec{\sigma}_h$, and therefore the total number of unknowns can be significantly reduced when compared to a classical LDG method. The coupling between elements is achieved solely by the hybrid variable $\vec{\lambda}_h$, and therefore, the system to be solved for globally can be condensed [10]. This means that the resulting system is usually much smaller than it would be for the standard LDG approach, which typically requires solving simple local problems on each element in an element-wise fashion.*

## 7. Numerical results: 2D examples

In this section we show two-dimensional numerical results. Here, we solve the (nonlinear) system of equations using Newton's method. The resulting linear system is solved using GMRES with block Jacobi preconditioning until the relative residual drops below $10^{-12}$. Newton's method is carried out until the $L_2$-norm of the residual drops below $10^{-10}$.

## 7.1. Linear advection equation

We first examine a system of linear advection equations. It can be written as in Eqn. (28). We compute a solution on $\Omega = [0, 2]^2$ at final time $T = 0.1$. The flux is chosen to be $\vec{f}(\vec{w}) = (\vec{f_1}, \vec{f_2})$ with

$$\vec{f_1}(\vec{w}) = \mathbf{A_1}\vec{w}, \quad \vec{f_2}(\vec{w}) = \mathbf{A_2}\vec{w}, \tag{32}$$

The vector of unknowns is $\vec{w} = (w_1, w_2)^T$. The matrices for this linear system are given by

$$\mathbf{A_1} = \begin{pmatrix} \frac{1}{3} & \frac{8}{3} \\ \frac{16}{3} & -\frac{7}{3} \end{pmatrix}, \ \mathbf{A_2} = \begin{pmatrix} -\frac{7}{3} & -\frac{5}{3} \\ -\frac{10}{3} & -\frac{2}{3} \end{pmatrix}. \tag{33}$$

These matrices have the same eigenvector basis, which means we can express these as $\mathbf{A_1} = \mathbf{S}\mathbf{D}_{A_1}\mathbf{S}^{-1}$ and $\mathbf{A_2} = \mathbf{S}\mathbf{D}_{A_2}\mathbf{S}^{-1}$ with

$$\mathbf{D}_{A_1} = \begin{pmatrix} -5 & 0 \\ 0 & 3 \end{pmatrix}, \mathbf{D}_{A_2} = \begin{pmatrix} 1 & 0 \\ 0 & -4 \end{pmatrix}, \mathbf{S} = \begin{pmatrix} -\frac{1}{2} & 1 \\ 1 & 1 \end{pmatrix}, \mathbf{S}^{-1} = \begin{pmatrix} -\frac{2}{3} & \frac{2}{3} \\ \frac{2}{3} & \frac{1}{3} \end{pmatrix}. \tag{34}$$

After choosing the initial conditions to be

$$\vec{w}_0(\vec{x}) = \begin{pmatrix} \sin(\pi(x+y)) \\ \sin(\pi(x+y)) \end{pmatrix}, \tag{35}$$

and taking into account periodic boundary conditions, we write the exact solution as

$$\vec{w}(\vec{x}, t) = \begin{pmatrix} \sin(\pi(x+y+t)) \\ \sin(\pi(x+y+t)) \end{pmatrix}. \tag{36}$$

We compute solutions on meshes that are presented in Fig. 8. Results are presented in Fig. 9 for the ratio $\frac{\Delta t}{\Delta x} = 0.025$. The errors for $w_1$ and $w_2$ are perfectly identical. The third order integrator reaches the expected order of convergence in all cases. For $p = 3$, the method is still third-order accurate, but it has a lower error than in the case with $p = 2$. The fourth-order integrator, however, does not achieve fourth-order in time. In the case $p < 3$, the method gets close to the expected order of $p + 1$ while for $p = 3$ the order deteriorates during the refinements. After the sixth refinement it seems not to converge any further. Most likely, this is behavior is observed due to stability issues of the fourth-order integrator.

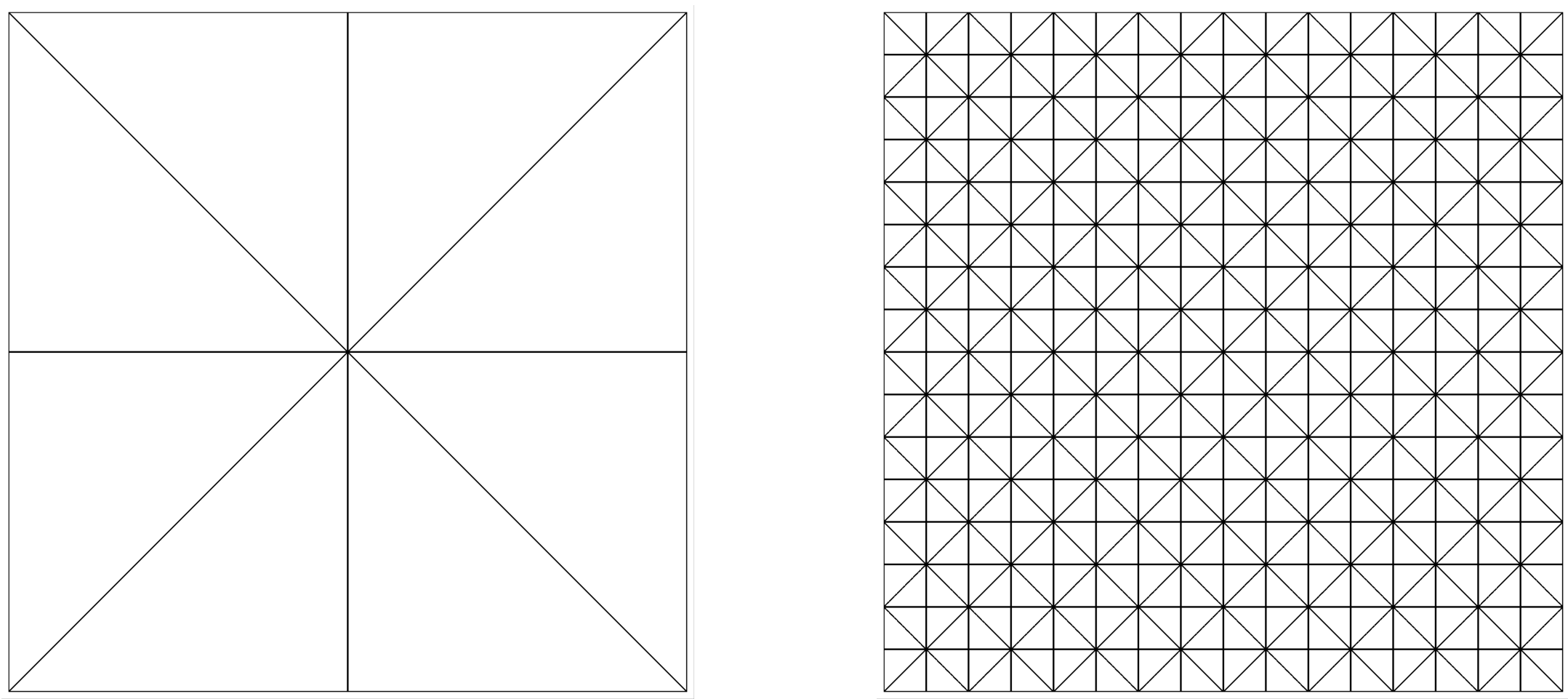

Figure 8: Left : Initial mesh. Right: Mesh after a total of three refinements.

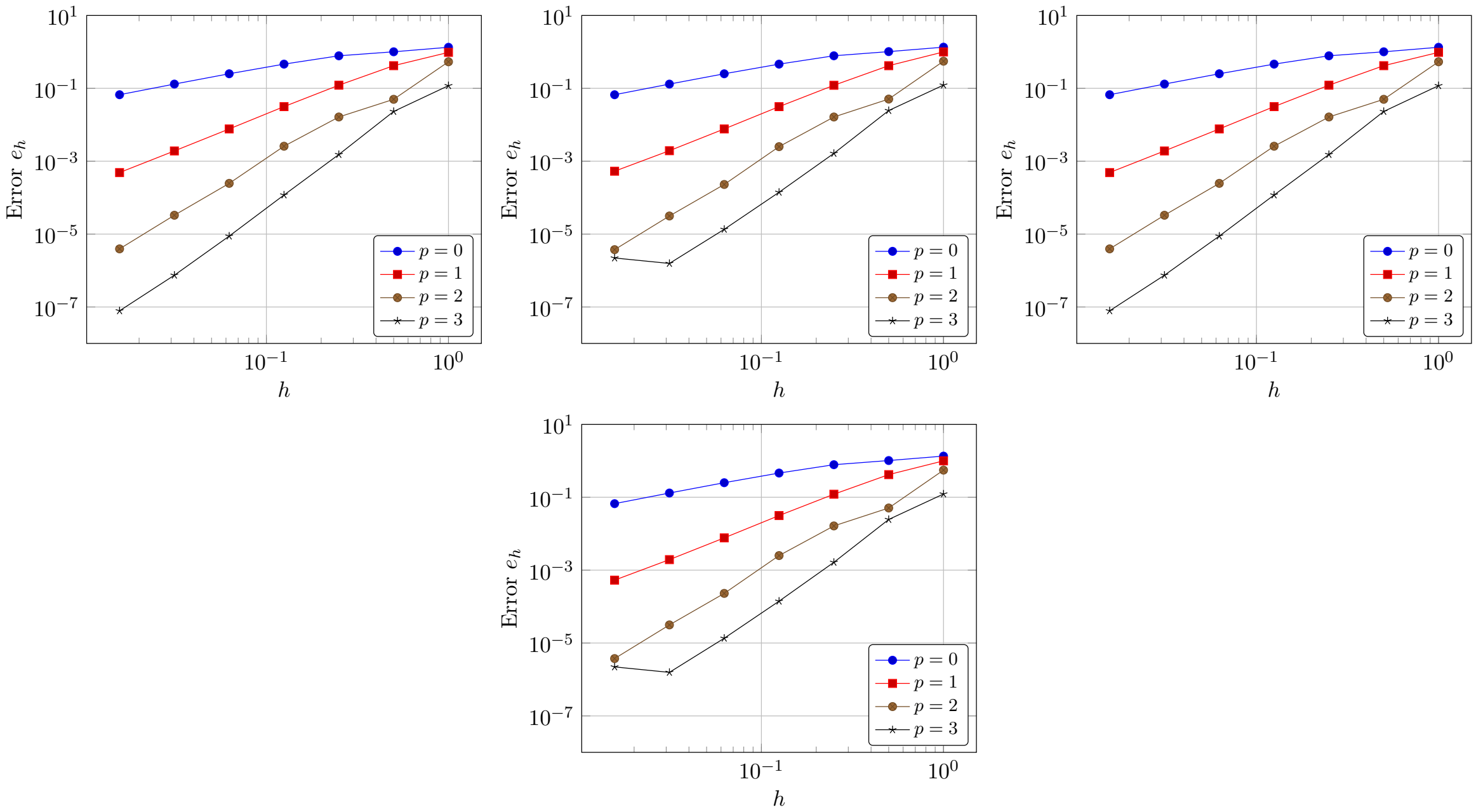

Figure 9: Numerical results for the linear coupled advection equation. Temporal integration is performed with the third-order integrator (7) (left) and the fourth-order (8) integrator (right). We show the results for components $w_1$ (top) and $w_2$ (bottom). In all computations we choose the ratio $\frac{\Delta t}{\Delta x}$ to be 0.025 to ensure stability of the numerical method. Tabulated results are given in Tab. 6 and 7.

### 7.1.1. Euler equations

As second test case in two space dimensions, we solve the Euler equations with periodic boundary conditions. The flux $\vec{f}(\vec{w}) = (\vec{f}_1, \vec{f}_2)$ is nonlinear and is given by

$$\begin{aligned} \vec{f}_1(\vec{w}) &= \left(\rho u, P + \rho u^2, \rho uv, u(E+P)\right)^T, \\ \vec{f}_2(\vec{w}) &= \left(\rho u, \rho uv, P + \rho v^2, v(E+P)\right)^T, \end{aligned} \tag{37}$$

and the vector of unknowns is $\vec{w} = (\rho, \rho u, \rho v, E)$, which define the density $\rho$, momentum $\rho u$ and $\rho v$ in the $x$- and $y$-direction, and energy $E$. The pressure $P$ is given by the equation of state

$$P = (\gamma - 1)\left(E - \frac{1}{2}\rho(u^2 + v^2)\right),$$

and the ratio of specific heats is $\gamma = 1.4$ for the test cases considered in this work. To analyze the accuracy of the method, we make use of a test case similar to the one presented in [31] that has a smooth analytical solution. The domain $\Omega = [0,2]^2$ is equipped with periodic boundary conditions, and the initial conditions are

$$\rho(x,y) = 1.0 + 0.2\sin(\pi(x+y)), \quad u = 0.7, \quad v = 0.3, \quad P = 1. \tag{38}$$

A convergence study is presented in Fig. 10, where we compute the solution to a final time of $T = 0.5$.

Both integrators produce very similar errors, but the third-order integrator has slightly lower errors. The conclusion is that the higher-order integrator does not exhibit any serious advantage over the lower-order integrator for this test case. For this problem, we find that increasing the polynomial order always increases the rate of convergence, which is in contrast to the previous cases. For example, in the linear test case, going from $p = 2$ to $p = 3$ decreased the error level, whereas the slope of the error graph stayed almost constant (cf. Fig. 9). For this problem, this actually *increases* the slope. Nevertheless, both integrators have a slight loss of convergence rate during refinements.

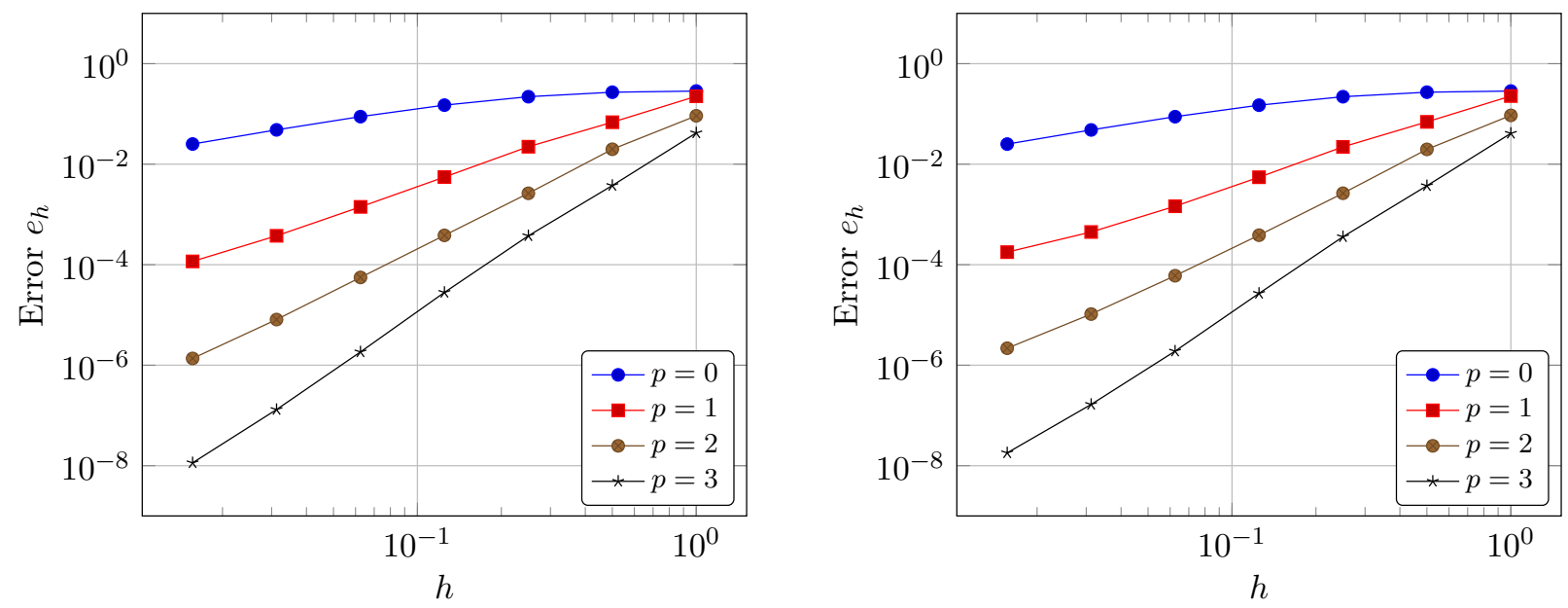

Figure 10: Numerical results for the Euler equations. Temporal integration is done via the third-order integrator (7) (left) and the fourth-order (8) integrator (right). We show the error in the density $\rho$. In all computations we choose the ratio $\frac{\Delta t}{\Delta x}$ to be 0.05 to ensure stability of the methods. Tabulated results are given in Tab. 8.

## 8. Conclusions and outlook

In this work, we present a novel application of high-order, implicit multiderivative time integrators to the discontinuous Galerkin framework. Two dimensional results are realized by employing the hybridized discontinuous Galerkin method in order to reduce the total number of unknowns that would otherwise be required to discretize the system. Results for a third- and fourth-order time integrator are presented, where we observe the expected order of convergence in time for all of our 1D test cases. For diffusion dominated problems, the integrators yield a stable scheme and the runtimes of the multiderivative schemes are reasonable compared to classical DIRK methods without the code being optimized for efficiency. However, despite the fact that we observe that the time integrators are not uniformly stable for convection dominated problems, these new methods work well for diffusion dominated problems and have reasonable runtimes. One possible explanation for this drawback is that higher derivatives carry negative coefficients, and effectively introduce anti-diffusion into the method. Future work must deal with improving stability properties of these methods. Furthermore, the extension of this methodology to the full Navier-Stokes equations - also for HDG - is of extreme importance and the subject of future work.

**Acknowledgements.** We would like to thank the anonymous reviewers for their helpful comments and suggestions to improve the quality of this work.

## A. Tabulated results

In this section we present complete error tables for the data presented summarily throughout the paper.

| $h$ | $p=0$ | order | $p=1$ | order | $p=2$ | order | $p=3$ | order |
|---|---|---|---|---|---|---|---|---|
| 5.00$e$−1 | 2.901$e$−1 | | 7.262$e$−2 | | 1.924$e$−2 | | 1.711$e$−2 | |
| 2.50$e$−1 | 7.445$e$−2 | 1.96 | 1.560$e$−2 | 2.22 | 2.838$e$−3 | 2.76 | 2.148$e$−3 | 2.99 |
| 1.25$e$−1 | 2.573$e$−2 | 1.53 | 3.725$e$−3 | 2.07 | 3.706$e$−4 | 2.94 | 2.892$e$−4 | 2.89 |
| 6.25$e$−2 | 1.154$e$−2 | 1.16 | 9.243$e$−4 | 2.01 | 4.786$e$−5 | 2.95 | 3.806$e$−5 | 2.93 |
| 3.13$e$−2 | 5.616$e$−3 | 1.04 | 2.306$e$−4 | 2.00 | 6.096$e$−6 | 2.97 | 4.899$e$−6 | 2.96 |
| 1.56$e$−2 | 2.790$e$−3 | 1.01 | 5.763$e$−5 | 2.00 | 7.697$e$−7 | 2.99 | 6.219$e$−7 | 2.98 |
| 7.81$e$−3 | 1.393$e$−3 | 1.00 | 1.440$e$−5 | 2.00 | 9.672$e$−8 | 2.99 | 7.837$e$−8 | 2.99 |
| 3.91$e$−3 | 6.960$e$−4 | 1.00 | 3.601$e$−6 | 2.00 | 1.212$e$−8 | 3.00 | 9.836$e$−9 | 2.99 |
| 1.95$e$−3 | 3.480$e$−4 | 1.00 | 9.002$e$−7 | 2.00 | 1.517$e$−9 | 3.00 | 1.232$e$−9 | 3.00 |
| 9.77$e$−4 | 1.740$e$−4 | 1.00 | 2.251$e$−7 | 2.00 | 1.898$e$−10 | 3.00 | 1.541$e$−10 | 3.00 |

| $h$ | $p=0$ | order | $p=1$ | order | $p=2$ | order | $p=3$ | order |
|---|---|---|---|---|---|---|---|---|
| 5.00$e$−1 | 3.384$e$−1 | | 6.631$e$−2 | | 1.879$e$−2 | | 6.051$e$−3 | |
| 2.50$e$−1 | 8.793$e$−2 | 1.94 | 1.631$e$−2 | 2.02 | 1.869$e$−3 | 3.33 | 3.218$e$−4 | 4.23 |
| 1.25$e$−1 | 2.744$e$−2 | 1.68 | 3.718$e$−3 | 2.13 | 2.314$e$−4 | 3.01 | 1.965$e$−5 | 4.03 |
| 6.25$e$−2 | 1.172$e$−2 | 1.23 | 9.221$e$−4 | 2.01 | 2.899$e$−5 | 3.00 | 1.222$e$−6 | 4.01 |
| 3.13$e$−2 | 5.636$e$−3 | 1.06 | 2.304$e$−4 | 2.00 | 3.627$e$−6 | 3.00 | 7.625$e$−8 | 4.00 |
| 1.56$e$−2 | 2.792$e$−3 | 1.01 | 5.761$e$−5 | 2.00 | 4.535$e$−7 | 3.00 | 4.764$e$−9 | 4.00 |
| 7.81$e$−3 | 1.393$e$−3 | 1.00 | 1.440$e$−5 | 2.00 | 5.669$e$−8 | 3.00 | 2.977$e$−10 | 4.00 |
| 3.91$e$−3 | 6.961$e$−4 | 1.00 | 3.601$e$−6 | 2.00 | 7.086$e$−9 | 3.00 | 1.862$e$−11 | 4.00 |
| 1.95$e$−3 | 3.480$e$−4 | 1.00 | 9.002$e$−7 | 2.00 | 8.858$e$−10 | 3.00 | 1.178$e$−12 | 3.98 |
| 9.77$e$−4 | 1.740$e$−4 | 1.00 | 2.251$e$−7 | 2.00 | 1.107$e$−10 | 3.00 | 1.286$e$−13 | 3.20 |

Table 1: Numerical results for the heat equation. Tabulated results of third-order (top) and forth-order (bottom) integrator as displayed in Fig. 2.

| $h$ | $p=0$ | order | $p=1$ | order | $p=2$ | order | $p=3$ | order |
|---|---|---|---|---|---|---|---|---|
| 5.00$e$−1 | 7.394$e$−1 | | 3.069$e$−1 | | 3.947$e$−1 | | 3.560$e$−1 | |
| 2.50$e$−1 | 6.846$e$−1 | 0.11 | 3.483$e$−1 | −0.18 | 9.526$e$−2 | 2.05 | 9.189$e$−2 | 1.95 |
| 1.25$e$−1 | 5.054$e$−1 | 0.44 | 1.812$e$−1 | 0.94 | 1.441$e$−2 | 2.72 | 1.411$e$−2 | 2.70 |
| 6.25$e$−2 | 3.283$e$−1 | 0.62 | 7.633$e$−2 | 1.25 | 1.875$e$−3 | 2.94 | 1.844$e$−3 | 2.94 |
| 3.13$e$−2 | 1.905$e$−1 | 0.79 | 2.837$e$−2 | 1.43 | 2.364$e$−4 | 2.99 | 2.328$e$−4 | 2.99 |
| 1.56$e$−2 | 1.028$e$−1 | 0.89 | 9.188$e$−3 | 1.63 | 2.961$e$−5 | 3.00 | 2.917$e$−5 | 3.00 |
| 7.81$e$−3 | 5.339$e$−2 | 0.95 | 2.655$e$−3 | 1.79 | 3.703$e$−6 | 3.00 | 3.649$e$−6 | 3.00 |
| 3.91$e$−3 | 2.720$e$−2 | 0.97 | 7.157$e$−4 | 1.89 | 4.629$e$−7 | 3.00 | 4.561$e$−7 | 3.00 |
| 1.95$e$−3 | 1.373$e$−2 | 0.99 | 1.859$e$−4 | 1.94 | 5.787$e$−8 | 3.00 | 5.702$e$−8 | 3.00 |
| 9.77$e$−4 | 6.897$e$−3 | 0.99 | 4.738$e$−5 | 1.97 | 7.234$e$−9 | 3.00 | 7.127$e$−9 | 3.00 |

| $h$ | $p=0$ | order | $p=1$ | order | $p=2$ | order | $p=3$ | order |
|---|---|---|---|---|---|---|---|---|
| 5.00$e$−1 | 6.988$e$−1 | | 5.194$e$−1 | | 7.035$e$−2 | | 2.520$e$−2 | |
| 2.50$e$−1 | 6.792$e$−1 | 0.04 | 1.243$e$−1 | 2.06 | 1.286$e$−2 | 2.45 | 1.282$e$−3 | 4.30 |
| 1.25$e$−1 | 5.152$e$−1 | 0.40 | 2.831$e$−2 | 2.13 | 1.664$e$−3 | 2.95 | 7.976$e$−5 | 4.01 |
| 6.25$e$−2 | 3.328$e$−1 | 0.63 | 6.771$e$−3 | 2.06 | 2.087$e$−4 | 3.00 | 5.032$e$−6 | 3.99 |
| 3.13$e$−2 | 1.915$e$−1 | 0.80 | 1.668$e$−3 | 2.02 | 2.611$e$−5 | 3.00 | 3.153$e$−7 | 4.00 |
| 1.56$e$−2 | 1.030$e$−1 | 0.89 | 4.154$e$−4 | 2.01 | 3.265$e$−6 | 3.00 | 1.970$e$−8 | 4.00 |
| 7.81$e$−3 | 5.341$e$−2 | 0.95 | 1.037$e$−4 | 2.00 | 4.081$e$−7 | 3.00 | 1.232$e$−9 | 4.00 |
| 3.91$e$−3 | 2.721$e$−2 | 0.97 | 2.593$e$−5 | 2.00 | 5.101$e$−8 | 3.00 | 7.696$e$−11 | 4.00 |
| 1.95$e$−3 | 1.373$e$−2 | 0.99 | 6.481$e$−6 | 2.00 | 6.377$e$−9 | 3.00 | 4.980$e$−12 | 3.95 |
| 9.77$e$−4 | 6.897$e$−3 | 0.99 | 1.620$e$−6 | 2.00 | 7.971$e$−10 | 3.00 | 3.961$e$−12 | 0.33 |

Table 2: Numerical results for the convection equation. Tabulated results of third-order (top) and fourth-order (bottom) integrator as displayed in Fig. 3.

| $h$ | $p=0$ | order | $p=1$ | order | $p=2$ | order | $p=3$ | order |
|---|---|---|---|---|---|---|---|---|
| $5.00e{-1}$ | $2.014e{-1}$ | | $8.409e{-2}$ | | $8.767e{-2}$ | | $1.133e{-1}$ | |
| $2.50e{-1}$ | $9.267e{-2}$ | 1.12 | $2.221e{-2}$ | 1.92 | $2.272e{-2}$ | 1.95 | $2.298e{-2}$ | 2.30 |
| $1.25e{-1}$ | $5.957e{-2}$ | 0.64 | $5.117e{-3}$ | 2.12 | $3.457e{-3}$ | 2.72 | $3.459e{-3}$ | 2.73 |
| $6.25e{-2}$ | $3.993e{-2}$ | 0.58 | $1.006e{-3}$ | 2.35 | $4.703e{-4}$ | 2.88 | $4.697e{-4}$ | 2.88 |
| $3.13e{-2}$ | $2.438e{-2}$ | 0.71 | $2.353e{-4}$ | 2.10 | $6.107e{-5}$ | 2.94 | $6.097e{-5}$ | 2.95 |
| $1.56e{-2}$ | $1.365e{-2}$ | 0.84 | $5.791e{-5}$ | 2.02 | $7.770e{-6}$ | 2.97 | $7.757e{-6}$ | 2.97 |
| $7.81e{-3}$ | $7.243e{-3}$ | 0.91 | $1.442e{-5}$ | 2.01 | $9.796e{-7}$ | 2.99 | $9.780e{-7}$ | 2.99 |
| $3.91e{-3}$ | $3.733e{-3}$ | 0.96 | $3.602e{-6}$ | 2.00 | $1.230e{-7}$ | 2.99 | $1.228e{-7}$ | 2.99 |
| $1.95e{-3}$ | $1.896e{-3}$ | 0.98 | $9.003e{-7}$ | 2.00 | $1.540e{-8}$ | 3.00 | $1.538e{-8}$ | 3.00 |
| $9.77e{-4}$ | $9.551e{-4}$ | 0.99 | $2.251e{-7}$ | 2.00 | $1.927e{-9}$ | 3.00 | $1.924e{-9}$ | 3.00 |

| $h$ | $p=0$ | order | $p=1$ | order | $p=2$ | order | $p=3$ | order |
|---|---|---|---|---|---|---|---|---|
| $5.00e{-1}$ | $1.332e{-1}$ | | $4.612e{-1}$ | | $6.349e{-2}$ | | $8.145e{-2}$ | |
| $2.50e{-1}$ | $1.039e{-1}$ | 0.36 | $7.866e{-2}$ | 2.55 | $5.563e{-3}$ | 3.51 | $5.594e{-3}$ | 3.86 |
| $1.25e{-1}$ | $7.032e{-2}$ | 0.56 | $6.107e{-3}$ | 3.69 | $4.147e{-4}$ | 3.75 | $3.671e{-4}$ | 3.93 |
| $6.25e{-2}$ | $4.524e{-2}$ | 0.64 | $9.880e{-4}$ | 2.63 | $3.642e{-5}$ | 3.51 | $2.330e{-5}$ | 3.98 |
| $3.13e{-2}$ | $2.620e{-2}$ | 0.79 | $2.333e{-4}$ | 2.08 | $3.874e{-6}$ | 3.23 | $1.462e{-6}$ | 3.99 |
| $1.56e{-2}$ | $1.418e{-2}$ | 0.89 | $5.786e{-5}$ | 2.01 | $4.608e{-7}$ | 3.07 | $9.149e{-8}$ | 4.00 |
| $7.81e{-3}$ | $7.385e{-3}$ | 0.94 | $1.443e{-5}$ | 2.00 | $5.687e{-8}$ | 3.02 | $5.719e{-9}$ | 4.00 |
| $3.91e{-3}$ | $3.770e{-3}$ | 0.97 | $3.604e{-6}$ | 2.00 | $7.089e{-9}$ | 3.00 | $3.575e{-10}$ | 4.00 |
| $1.95e{-3}$ | $1.905e{-3}$ | 0.98 | $9.007e{-7}$ | 2.00 | $8.857e{-10}$ | 3.00 | $2.235e{-11}$ | 4.00 |
| $9.77e{-4}$ | $9.575e{-4}$ | 0.99 | $2.251e{-7}$ | 2.00 | $1.107e{-10}$ | 3.00 | $1.416e{-12}$ | 3.98 |

Table 3: Numerical results for the convection-diffusion equation with smooth initial condition. Tabulated results of third-order (top) and fourth-order (bottom) integrator as displayed in Fig. 4.

| $h$ | error | order | $time[s]$ |
|---|---|---|---|
| 5.000$e$−1 | 1.133$e$−1 | | 1.371$e$−2 |
| 2.500$e$−1 | 2.298$e$−2 | 2.30 | 2.054$e$−2 |
| 1.250$e$−1 | 3.459$e$−3 | 2.73 | 4.390$e$−2 |
| 6.250$e$−2 | 4.697$e$−4 | 2.88 | 1.011$e$−1 |
| 3.125$e$−2 | 6.097$e$−5 | 2.95 | 2.526$e$−1 |
| 1.563$e$−2 | 7.757$e$−6 | 2.97 | 7.274$e$−1 |
| 7.813$e$−3 | 9.780$e$−7 | 2.99 | 2.321$e$+0 |
| 3.906$e$−3 | 1.228$e$−7 | 2.99 | 8.253$e$+0 |
| 1.953$e$−3 | 1.538$e$−8 | 3.00 | 3.186$e$+1 |
| 9.766$e$−4 | 1.924$e$−9 | 3.00 | 1.286$e$+2 |

| $h$ | error | order | $time[s]$ |
|---|---|---|---|
| 5.000$e$−1 | 1.429$e$−1 | | 6.308$e$−3 |
| 2.500$e$−1 | 3.226$e$−2 | 2.15 | 9.798$e$−3 |
| 1.250$e$−1 | 5.426$e$−3 | 2.57 | 2.268$e$−2 |
| 6.250$e$−2 | 8.008$e$−4 | 2.76 | 5.643$e$−2 |
| 3.125$e$−2 | 1.088$e$−4 | 2.88 | 1.635$e$−1 |
| 1.563$e$−2 | 1.416$e$−5 | 2.94 | 5.522$e$−1 |
| 7.813$e$−3 | 1.805$e$−6 | 2.97 | 2.025$e$+0 |
| 3.906$e$−3 | 2.277$e$−7 | 2.99 | 7.777$e$+0 |
| 1.953$e$−3 | 2.860$e$−8 | 2.99 | 3.193$e$+1 |
| 9.766$e$−4 | 3.583$e$−9 | 3.00 | 1.342$e$+2 |

| $h$ | error | order | $time[s]$ |
|---|---|---|---|
| 5.000$e$−1 | 8.145$e$−2 | | 1.032$e$−2 |
| 2.500$e$−1 | 5.594$e$−3 | 3.86 | 2.017$e$−2 |
| 1.250$e$−1 | 3.671$e$−4 | 3.93 | 4.359$e$−2 |
| 6.250$e$−2 | 2.330$e$−5 | 3.98 | 9.855$e$−2 |
| 3.125$e$−2 | 1.462$e$−6 | 3.99 | 2.461$e$−1 |
| 1.563$e$−2 | 9.149$e$−8 | 4.00 | 6.854$e$−1 |
| 7.813$e$−3 | 5.719$e$−9 | 4.00 | 2.206$e$+0 |
| 3.906$e$−3 | 3.575$e$−10 | 4.00 | 7.756$e$+0 |
| 1.953$e$−3 | 2.235$e$−11 | 4.00 | 2.999$e$+1 |
| 9.766$e$−4 | 1.416$e$−12 | 3.98 | 1.201$e$+2 |

| $h$ | error | order | $time[s]$ |
|---|---|---|---|
| 5.000$e$−1 | 4.906$e$−2 | | 5.399$e$−3 |
| 2.500$e$−1 | 3.520$e$−3 | 3.80 | 1.162$e$−2 |
| 1.250$e$−1 | 2.282$e$−4 | 3.95 | 2.720$e$−2 |
| 6.250$e$−2 | 1.434$e$−5 | 3.99 | 7.142$e$−2 |
| 3.125$e$−2 | 8.957$e$−7 | 4.00 | 2.167$e$−1 |
| 1.563$e$−2 | 5.591$e$−8 | 4.00 | 7.741$e$−1 |
| 7.813$e$−3 | 3.491$e$−9 | 4.00 | 2.915$e$+0 |
| 3.906$e$−3 | 2.181$e$−10 | 4.00 | 1.156$e$+1 |
| 1.953$e$−3 | 1.364$e$−11 | 4.00 | 4.757$e$+1 |
| 9.766$e$−4 | 8.726$e$−13 | 3.97 | 2.091$e$+2 |

Table 4: Results for the convection-diffusion equation. Tables with errors and computing times used in Fig. 5. In clockwise direction starting at the top left the results of the third order integrator, the DIRK method of Cash, the DIRK method of Hairer and Wanner and the fourth order integrator are shown.

| $h$ | $p=0$ | order | $p=1$ | order | $p=2$ | order | $p=3$ | order |
|---|---|---|---|---|---|---|---|---|
| 5.00$e$−1 | 1.254$e$−2 | | 1.484$e$−1 | | 1.900$e$−2 | | 2.449$e$−2 | |
| 2.50$e$−1 | 7.301$e$−2 | −2.54 | 1.726$e$−2 | 3.10 | 2.867$e$−3 | 2.73 | 2.269$e$−3 | 3.43 |
| 1.25$e$−1 | 5.596$e$−2 | 0.38 | 3.754$e$−3 | 2.20 | 3.197$e$−4 | 3.17 | 2.198$e$−4 | 3.37 |
| 6.25$e$−2 | 3.761$e$−2 | 0.57 | 8.954$e$−4 | 2.07 | 3.702$e$−5 | 3.11 | 2.285$e$−5 | 3.27 |
| 3.13$e$−2 | 2.239$e$−2 | 0.75 | 2.203$e$−4 | 2.02 | 4.300$e$−6 | 3.11 | 2.280$e$−6 | 3.33 |
| 1.56$e$−2 | 1.230$e$−2 | 0.86 | 5.481$e$−5 | 2.01 | 5.128$e$−7 | 3.07 | 2.345$e$−7 | 3.28 |
| 7.81$e$−3 | 6.454$e$−3 | 0.93 | 1.369$e$−5 | 2.00 | 6.251$e$−8 | 3.04 | 2.564$e$−8 | 3.19 |
| 3.91$e$−3 | 3.308$e$−3 | 0.96 | 3.421$e$−6 | 2.00 | 7.717$e$−9 | 3.02 | 2.959$e$−9 | 3.12 |
| 1.95$e$−3 | 1.675$e$−3 | 0.98 | 8.551$e$−7 | 2.00 | 9.587$e$−10 | 3.01 | 3.541$e$−10 | 3.06 |

| $h$ | $p=0$ | order | $p=1$ | order | $p=2$ | order | $p=3$ | order |
|---|---|---|---|---|---|---|---|---|
| 5.00$e$−1 | 8.023$e$−2 | | 2.948$e$−1 | | 6.399$e$−2 | | 7.559$e$−2 | |
| 2.50$e$−1 | 8.767$e$−2 | −0.13 | 1.677$e$−2 | 4.14 | 4.270$e$−3 | 3.91 | 3.981$e$−3 | 4.25 |
| 1.25$e$−1 | 6.395$e$−2 | 0.46 | 3.247$e$−3 | 2.37 | 2.816$e$−4 | 3.92 | 1.858$e$−4 | 4.42 |
| 6.25$e$−2 | 4.063$e$−2 | 0.65 | 8.045$e$−4 | 2.01 | 2.972$e$−5 | 3.24 | 1.241$e$−5 | 3.90 |
| 3.13$e$−2 | 2.334$e$−2 | 0.80 | 2.072$e$−4 | 1.96 | 3.583$e$−6 | 3.05 | 8.049$e$−7 | 3.95 |
| 1.56$e$−2 | 1.256$e$−2 | 0.89 | 5.307$e$−5 | 1.96 | 4.486$e$−7 | 3.00 | 5.090$e$−8 | 3.98 |
| 7.81$e$−3 | 6.525$e$−3 | 0.94 | 1.346$e$−5 | 1.98 | 5.644$e$−8 | 2.99 | 3.193$e$−9 | 3.99 |
| 3.91$e$−3 | 3.326$e$−3 | 0.97 | 3.392$e$−6 | 1.99 | 7.088$e$−9 | 2.99 | 1.997$e$−10 | 4.00 |
| 1.95$e$−3 | 1.679$e$−3 | 0.99 | 8.515$e$−7 | 1.99 | 8.883$e$−10 | 3.00 | 1.248$e$−11 | 4.00 |

Table 5: Numerical results for Burgers equation. Tabulated results of third-order (top) and fourth-order (bottom) integrator as displayed in Fig. 7.

| $h$ | $p=0$ | order | $p=1$ | order | $p=2$ | order | $p=3$ | order |
|---|---|---|---|---|---|---|---|---|
| 1.00e+0 | 1.338e+0 | 0.00 | 9.699e−1 | 0.00 | 5.361e−1 | | 1.184e−1 | |
| 5.00e−1 | 1.007e+0 | 0.41 | 4.185e−1 | 0.41 | 4.997e−2 | 3.42 | 2.320e−2 | 2.35 |
| 2.50e−1 | 7.812e−1 | 0.37 | 1.219e−1 | 0.37 | 1.636e−2 | 1.61 | 1.529e−3 | 3.92 |
| 1.25e−1 | 4.619e−1 | 0.76 | 3.152e−2 | 0.76 | 2.579e−3 | 2.67 | 1.184e−4 | 3.69 |
| 6.25e−2 | 2.511e−1 | 0.88 | 7.681e−3 | 0.88 | 2.467e−4 | 3.39 | 8.813e−6 | 3.75 |
| 3.13e−2 | 1.309e−1 | 0.94 | 1.915e−3 | 0.94 | 3.295e−5 | 2.90 | 7.451e−7 | 3.56 |
| 1.56e−2 | 6.680e−2 | 0.97 | 4.885e−4 | 0.97 | 3.932e−6 | 3.07 | 7.863e−8 | 3.24 |

| $h$ | $p=0$ | order | $p=1$ | order | $p=2$ | order | $p=3$ | order |
|---|---|---|---|---|---|---|---|---|
| 1.00e+0 | 1.338e+0 | 0.00 | 9.699e−1 | 0.00 | 5.361e−1 | | 1.184e−1 | |
| 5.00e−1 | 1.007e+0 | 0.41 | 4.185e−1 | 0.41 | 4.997e−2 | 3.42 | 2.320e−2 | 2.35 |
| 2.50e−1 | 7.812e−1 | 0.37 | 1.219e−1 | 0.37 | 1.636e−2 | 1.61 | 1.529e−3 | 3.92 |
| 1.25e−1 | 4.619e−1 | 0.76 | 3.152e−2 | 0.76 | 2.579e−3 | 2.67 | 1.184e−4 | 3.69 |
| 6.25e−2 | 2.511e−1 | 0.88 | 7.681e−3 | 0.88 | 2.467e−4 | 3.39 | 8.813e−6 | 3.75 |
| 3.13e−2 | 1.309e−1 | 0.94 | 1.915e−3 | 0.94 | 3.295e−5 | 2.90 | 7.451e−7 | 3.56 |
| 1.56e−2 | 6.680e−2 | 0.97 | 4.885e−4 | 0.97 | 3.932e−6 | 3.07 | 7.863e−8 | 3.24 |

Table 6: Numerical results for the linear coupled advection equation. Tabulated results of first (top) and second (bottom) component obtained by third-order integrator displayed in Fig. 9.

| $h$ | $p=0$ | order | $p=1$ | order | $p=2$ | order | $p=3$ | order |
|---|---|---|---|---|---|---|---|---|
| 1.00e+0 | 1.347e+0 | 0.00 | 9.976e−1 | 0.00 | 5.557e−1 | | 1.229e−1 | |
| 5.00e−1 | 1.018e+0 | 0.40 | 4.180e−1 | 0.40 | 5.075e−2 | 3.45 | 2.458e−2 | 2.32 |
| 2.50e−1 | 7.820e−1 | 0.38 | 1.212e−1 | 0.38 | 1.635e−2 | 1.63 | 1.641e−3 | 3.90 |
| 1.25e−1 | 4.615e−1 | 0.76 | 3.146e−2 | 0.76 | 2.510e−3 | 2.70 | 1.407e−4 | 3.54 |
| 6.25e−2 | 2.509e−1 | 0.88 | 7.683e−3 | 0.88 | 2.301e−4 | 3.45 | 1.354e−5 | 3.38 |
| 3.13e−2 | 1.307e−1 | 0.94 | 1.944e−3 | 0.94 | 3.142e−5 | 2.87 | 1.566e−6 | 3.11 |
| 1.56e−2 | 6.674e−2 | 0.97 | 5.319e−4 | 0.97 | 3.745e−6 | 3.07 | 2.211e−6 | −0.50 |

| $h$ | $p=0$ | order | $p=1$ | order | $p=2$ | order | $p=3$ | order |
|---|---|---|---|---|---|---|---|---|
| 1.00e+0 | 1.347e+0 | 0.00 | 9.976e−1 | 0.00 | 5.557e−1 | | 1.229e−1 | |
| 5.00e−1 | 1.018e+0 | 0.40 | 4.180e−1 | 0.40 | 5.075e−2 | 3.45 | 2.458e−2 | 2.32 |
| 2.50e−1 | 7.820e−1 | 0.38 | 1.212e−1 | 0.38 | 1.635e−2 | 1.63 | 1.641e−3 | 3.90 |
| 1.25e−1 | 4.615e−1 | 0.76 | 3.146e−2 | 0.76 | 2.510e−3 | 2.70 | 1.407e−4 | 3.54 |
| 6.25e−2 | 2.509e−1 | 0.88 | 7.683e−3 | 0.88 | 2.301e−4 | 3.45 | 1.354e−5 | 3.38 |
| 3.13e−2 | 1.307e−1 | 0.94 | 1.944e−3 | 0.94 | 3.142e−5 | 2.87 | 1.566e−6 | 3.11 |
| 1.56e−2 | 6.674e−2 | 0.97 | 5.319e−4 | 0.97 | 3.745e−6 | 3.07 | 2.211e−6 | −0.50 |

Table 7: Numerical results for the linear coupled advection equation. Tabulated results of first (top) and second (bottom) component obtained by fourth-order integrator Fig. 9.

| $h$ | $p = 0$ | order | $p = 1$ | order | $p = 2$ | order | $p = 3$ | order |
|---|---|---|---|---|---|---|---|---|
| 1.00$e$+0 | 2.851$e$−1 | 0.00 | 2.242$e$−1 | 0.00 | 9.174$e$−2 | | 4.183$e$−2 | |
| 5.00$e$−1 | 2.704$e$−1 | 0.09 | 6.833$e$−2 | 0.08 | 1.982$e$−2 | 2.21 | 3.760$e$−3 | 3.48 |
| 2.50$e$−1 | 2.202$e$−1 | 0.32 | 2.226$e$−2 | 0.30 | 2.633$e$−3 | 2.91 | 3.757$e$−4 | 3.32 |
| 1.25$e$−1 | 1.495$e$−1 | 0.52 | 5.542$e$−3 | 0.56 | 3.833$e$−4 | 2.78 | 2.805$e$−5 | 3.74 |
| 6.25$e$−2 | 8.835$e$−2 | 0.75 | 1.413$e$−3 | 0.76 | 5.574$e$−5 | 2.78 | 1.836$e$−6 | 3.93 |
| 3.13$e$−2 | 4.817$e$−2 | 0.88 | 3.756$e$−4 | 0.88 | 8.096$e$−6 | 2.78 | 1.310$e$−7 | 3.81 |
| 1.56$e$−2 | 2.517$e$−2 | 0.94 | 1.162$e$−4 | 0.94 | 1.363$e$−6 | 2.57 | 1.141$e$−8 | 3.52 |

| $h$ | $p = 0$ | order | $p = 1$ | order | $p = 2$ | order | $p = 3$ | order |
|---|---|---|---|---|---|---|---|---|
| 1.00$e$+0 | 2.852$e$−1 | 0.00 | 2.269$e$−1 | 0.00 | 9.354$e$−2 | | 4.102$e$−2 | |
| 5.00$e$−1 | 2.707$e$−1 | 0.09 | 6.949$e$−2 | 0.08 | 1.971$e$−2 | 2.25 | 3.720$e$−3 | 3.46 |
| 2.50$e$−1 | 2.200$e$−1 | 0.32 | 2.211$e$−2 | 0.30 | 2.639$e$−3 | 2.90 | 3.593$e$−4 | 3.37 |
| 1.25$e$−1 | 1.492$e$−1 | 0.52 | 5.497$e$−3 | 0.56 | 3.863$e$−4 | 2.77 | 2.683$e$−5 | 3.74 |
| 6.25$e$−2 | 8.810$e$−2 | 0.76 | 1.454$e$−3 | 0.76 | 6.029$e$−5 | 2.68 | 1.902$e$−6 | 3.82 |
| 3.13$e$−2 | 4.803$e$−2 | 0.88 | 4.472$e$−4 | 0.88 | 1.040$e$−5 | 2.53 | 1.656$e$−7 | 3.52 |
| 1.56$e$−2 | 2.509$e$−2 | 0.94 | 1.777$e$−4 | 0.94 | 2.176$e$−6 | 2.26 | 1.800$e$−8 | 3.20 |

Table 8: Numerical results for the Euler equations. Tabulated results of third-order (top) and fourth-order (bottom) integrator as displayed in Fig. 10.